\documentclass{gtpart}     
\usepackage{graphicx}  
\usepackage[all]{xy}
\title{Minimal Entropy and Geometric Decompositions in Dimension Four}
\author[P. Su\'arez-Serrato \; \today]{P. Su\'arez-Serrato}

\givenname{Pablo}
\surname{Su\'arez-Serrato}
\address{Mathematisches Institut LMU, Theresienstrasse 39, M\"unchen 80333 Deutschland.}
\email{pablo@math.lmu.de}

\keyword{minimal entropy, geodesic flows, geometric structures}
\subject{primary}{msc2000}{37B40}
\subject{secondary}{msc2000}{53D25}

\arxivreference{math.GT/0601759}
\arxivpassword{pdxcm}

\volumenumber{}
\issuenumber{}
\publicationyear{}
\papernumber{}
\startpage{}
\endpage{}
\doi{}
\MR{}
\Zbl{}
\received{}
\revised{}
\accepted{}
\published{}
\publishedonline{}
\proposed{}
\seconded{}
\corresponding{}
\editor{}
\version{}

\newtheorem*{A}{Theorem A}

\theoremstyle{definition}
\newtheorem*{rem}{Remark}          

\newtheorem{Theorem}{Theorem}
\newtheorem{Lemma}[Theorem]{Lemma}
\newtheorem{Corollary}[Theorem]{Corollary}

\newtheorem{Proposition}[Theorem]{Proposition}
\newtheorem{Conjecture}[Theorem]{Conjecture}
\newtheorem{Definition}[Theorem]{Definition}

\def \dim{{\mbox {dim}}\,}

\def\Z{{\mathbb Z}}

\def \re{{\mathbb R}}

\def \C{{\mathbb C}}

\def \H{{\mathbb H}}

\def \0{\lambda_{0}}

\def\h{{\rm h}_{\rm top}(g)}

\def\Fs{{\mathcal F}}
\def\Ts{{\mathcal T}}
\def\Sl{\widetilde{\rm SL}_{2}}
\def\is{\rm{Iso}}
\def \X{{\mathbb X}}
\def \E{\mathbb{E}}

\begin{document}

\begin{abstract} We show vanishing results about the infimum of the topological entropy of the geodesic flow of homogeneous smooth four manifolds. We prove that any closed oriented \emph{geometric} four manifold has zero minimal entropy if and only if it has zero simplicial volume. We also show that if a four manifold $M$ admits a geometric decomposition, in the sense of Thurston, and does not have geometric pieces modelled on hyperbolic four-space $\H^4$, the complex hyperbolic plane $\H_{\C}^2$ or the product of two hyperbolic planes $\H^2 \times \H^2$,  then $M$ admits an $\Fs$-structure. It follows that $M$ has zero minimal entropy and collapses with curvature bounded from below. We then analyse whether or not $M$ admits a metric whose topological entropy coincides with the minimal entropy of $M$ and provide new examples of manifolds for which the minimal entropy problem cannot be solved.

 \end{abstract}

\maketitle


\section{Introduction}

A \emph{model geometry}, in the sense of W \,P Thurston,  is a complete simply connected Riemannian manifold $X$ such that the group of isometries acts transitively on $X$ and contains a discrete subgroup with a finite volume quotient. The maximal four dimensional geometries were classified by R Filipkiewicz \cite{Fi}. In this note we will focus on the minimal entropy problem for smooth 4-manifolds $M$ which are geometrisable in the sense of Thurston; $M$ is diffeomorphic to a connected sum of manifolds which admit a decomposition into pieces which are modelled on a Thurston geometry.

The minimal entropy ${\rm h}(M)$ of a closed smooth manifold $M$ is the infimum of the topological entropy $\h$ of the geodesic flow of $g$ over the familiy of $C^{\infty}$ Riemannian metrics on $M$ with unit volume. A metric $g_0$ is {\it entropy minimising} if it achieves this infimum ${\rm h_{top}}(g_0)= {\rm h}(M)$, when such a metric exists we say the minimal entropy problem can be solved for $M$.

The minimal entropy $ { \rm h}(M)$ of an $n$--manifold $M$ is related to its simplicial volume $||M||$, volume entropy $\lambda (M)$ and minimal volume $ {\rm MinVol}(M)$ according to the inequalities noticed by M Gromov \cite{Gro}, A Manning \cite{Man} and by various authors \cite{BCG,Pat},
$$\frac{n^{n/2}}{n!}||M|| \leq \lambda(M)^{n} \leq  {\rm h} (M)^{n} \leq (n-1)^{n} {\rm MinVol}(M). $$

The simplicial and minimal volumes were defined by Gromov in the seminal paper \cite{Gro}. Both the simplicial volume and volume entropy are known to be homotopy invariant, see I \, K Babenko \cite{Bab} and M Brunnbauer \cite{Br}. However, L Bessi\`{e}res \cite{Bes} has shown that the  minimal volume ${\rm MinVol}(M)$ depends on the differentiable structure of $M$. In fact D Kotschick has proven that even the vanishing of ${\rm MinVol}(M)$ depends on the differentiable structure \cite{K2}. As the question of whether the minimal entropy is a homotopy invariant is still unresolved, it is interesting to calculate it and compare it with the invariants mentioned above.

The instrument we will use to show that these invariants vanish is a generalisation of a local torus action,  called an $\mathcal{F}$--structure. J Cheeger and Gromov showed in \cite{CG} that if a manifold $M$ admits a \emph{polarised} $\mathcal{F}$--structure then ${\rm MinVol}(M)=0$. The simplest example of a polarised $\mathcal{F}$--structure is a free $S^1$--action on $M$. Paternain and J Petean  proved that if $M$ admits any $\mathcal{F}$--structure then ${\rm h} (M)=0$,  $M$ collapses with curvature bounded from below and the Yamabe invariant of $M$ is non-negative \cite{PP}.

In dimension four there exist smooth manifolds that admit $\Fs$--structures and which are homeomorphic to manifolds that do not admit them, see Paternain and Petean \cite{PP} and C LeBrun \cite{LB}. The results in this paper provide a basis of examples with which to compare manifolds in the same homeomorphism class.

The relevant definitions will be reviewed in the following sections.

Let ${\bf H}$ and ${\bf V}$ be the following sets of four dimensional geometries:
\[ { \bf H }=\{\mathbb{H}^{4}, \mathbb{H}^{2}\times\mathbb{H}^{2}, \mathbb{H}_{\C}^{2} \} \]
\[ {\bf V}  =  \left\{
\begin{array}{cccc}
     \mathbb{S}^{4} & \mathbb{C}{\rm P}^{2} & \mathbb{S}^{3}\times \mathbb{E} & \mathbb{H}^{3}\times \mathbb{E}   \\
 \Sl \times \mathbb{E} & \mathbb{N}il^{3}\times \mathbb{E}  & \mathbb{N}il^{4} &  \mathbb{S}ol^{4}_{1} \\
 \mathbb{S}^{2} \times\mathbb{E}^{2} & \mathbb{H}^{2}\times\mathbb{E}^{2}
& \mathbb{S}ol^{4}_{m,n} & \mathbb{S}ol^{4}_{0}  \\
 \mathbb{S}^{2}\times \mathbb{S}^{2}  & \mathbb{S}^{2}\times\mathbb{H}^{2} & \mathbb{E}^{4} & \mathbb{F}^4
\end{array} \right\}
\]

Together ${\bf H}$ and ${\bf V}$ constitute all the four--dimensional geometries that admit finite volume quotients.

Our first result is about compact manifolds which are modelled on a single geometry, these are called {\emph{geometric}} manifolds.

\medskip

{\bf Theorem~A.} {\it Let $M$ be a smooth oriented and closed geometric four manifold. The following notions are equivalent:

i) $M$ is modelled on a geometry in ${\bf V}$

ii) $M$ admits an $\Fs$--structure

iii) $M$ has zero minimal entropy, ${\rm h}(M)=0$

iv) The simplicial volume of $M$ vanishes, $||M||=0$

v)  $M$ collapses with curvature bounded from below,  ${\rm Vol_{K}}(M)=0$ }
\medskip

The main novel ingredient here is the proof that {\em i)} implies {\em ii)}, this is shown in Theorem 1. The other equivalences follow from the inequalities above and an application of results of Gromov and Thurston, as shown in Proposition \ref{Hpos}. This can be seen in the next diagram, which summarises the proof of Theorem A.

\[ \xymatrix{  & *\txt{ \sl i) $M$ is modelled\\ \sl on a geometry in ${\bf V}$}\ar@{=>}[ldd]_-{\txt{Theorem 1}} &  \\ & &\\
      *\txt{ \sl ii) $M$ admits\\ \sl an $\Fs$--structure }\ar@{=>}[dd]_-{\txt{Paternain--Petean\\Theorem}}  & & *\txt{ \sl v)  $M$ collapses with\\ \sl curvature bounded\\ \sl from below }\ar@{=>}[uul]_-{\txt{Contrapositive to\\Proposition 2}}  \\ & & \\
 *\txt{ \sl iii) ${\rm h}(M)=0$}\ar@{=>}[rr]_-{\txt{Inequalities of\\asymptotic invariants}} & &    *\txt{ \sl  iv) $||M||=0$ }\ar@{=>}[uu]_-{\txt{Theorem 1\\and Proposition 2}}} \]

%
%

The definitions and properties of these terms will appear in the following sections.


 Relying on results of  Thurston and Gromov \cite{Gro} we can then see that the contents of  Theorem A can be rephrased in the following way.
\\

\noindent {\it  Let $M$ be an oriented closed smooth geometric four manifold, the following are equivalent:

i) ${\rm h}(M)>0$

ii) ${\rm Vol_{K}} (M) > 0$

iii) $||M|| \neq 0$

iv) $M$ does not admit any $\Fs$--structures

v) $M$ is modelled on a geometry in ${\bf H}$}
\medskip

J\, W Anderson and Paternain showed  in \cite[Theorem 2.9]{AP} that  for a geometric 3--manifold $M$ it is equivalent for its simplicial volume, minimal entropy or minimal volume to vanish and for $M$ to be a {\emph{graph manifold}}. If a geometric 3--manifold $M$ admits a geometric structure modelled on a geometry which is not $\mathbb{H}^3$ then $M$ is a graph manifold. By the results  of G Besson, G Courtois and S Gallot \cite{BCG} if $M$ is modelled on $\mathbb{H}^3$ then the minimal entropy of $M$ is strictly positive and it is achieved by the hyperbolic metric.

In the same vein Theorem A shows that vanishing of the minimal entropy is an obstruction to the manifold being of \emph{hyperbolic type} in the extended sense of it being modelled on a geometry in ${\bf H}$.

A manifold $M$ is said to admit a \emph{geometric decomposition} if it admits a finite collection of 2--sided hypersurfaces $S$ such that each component of $M-S$ is geometric. A manifold is \emph{geometrisable \'a la Thurston} if it is diffeomorphic to a connected sum of manifolds with a geometric decomposition. After inspecting every possible geometric decomposition and every type of geometrizable smooth four manifold, we can extend Theorem A to the geometrisable case. The main result of this paper is the following theorem.

\medskip
\noindent {\bf Theorem~B.} {\it  Let $M$ be a closed orientable smooth four manifold which is geometrisable \`a la Thurston. If all of the geometric pieces of $M$ are modelled on geometries in ${\bf V}$ then $M$ admits an $\Fs$-structure. Consequently ${\rm h} (M) = {\rm Vol_{K}} (M) = 0$.}
\medskip

Therefore all closed geometrisable smooth four manifolds $M$ which are known to have $\| M \|=0$ also have ${\rm h} (M) =0$. It should be noted that there are no known examples of manifolds with zero simplicial volume and positive minimal entropy. So Theorem B also shows that such an example can not be constructed in dimension four by means of geometric decompositions.

The minimal entropy problem for geometric four manifolds has been treated by Paternain and Petean in \cite{PP1}. They have shown that if $M$ admits a geometric structure modelled on $\mathbb{S}^{4},
\mathbb{C}{\rm P}^{2}, \mathbb{S}^{3}\times \mathbb{E} ,
 \mathbb{N}il^{3}\times \mathbb{E} ,
\mathbb{N}il^{4}, \mathbb{S}^{2} \times
\mathbb{E}^{2},
 \mathbb{S}^{2}\times \mathbb{S}^{2}$ or  $\mathbb{E}^{4}$ then $M$ admits a metric with zero topological entropy. Whereas if $M$ is modelled on $ \mathbb{S}^{2}\times \mathbb{H}^{2},  \mathbb{H}^{3}\times \mathbb{E} , \Sl\times \mathbb{E} , \mathbb{H}^{2}\times\mathbb{E}^{2},  \mathbb{S}ol^{4}_{1}, \mathbb{S}ol^{4}_{0} $ or $\mathbb{S}ol^{4}_{m,n}$ then the fundamental group of $M$ has exponential growth, which implies that any  smooth metric on $M$ has positive topological entropy by a result of Manning \cite{Man}. It follows from Theorem A that the minimal entropy problem cannot be solved for a manifold $M$ modelled on any of these geometries, since we can endow them with $\mathcal{T}$--structures.

On the other hand, for manifolds modelled on $\mathbb{H}^4$ and $\mathbb{H}_{\C}^2$ the work of Besson, Courtois and Gallot implies that the minimal entropy problem is solved by their respective \emph{hyperbolic} and \emph{locally symmetric} metrics \cite{BCG}. Moreover, finite volume manifolds modelled on these two geometries have positive simplicial volume.

Solving the minimal entropy problem for manifolds $M$ modelled on the geometry $\mathbb{H}^2\times \mathbb{H}^2$ remains open. A possible candidate for an entropy minimising metric on $M$ could be the product metric on $\mathbb{H}^2\times \mathbb{H}^2$ inherited on $M$ as a quotient.

We show in Lemma \ref{d-eg} that if an orientable  $4$--manifold $M$ admits a proper geometric decomposition, its fundamental group $\pi_1(M)$ is not trivial and grows exponentially. So the only manifolds considered in Theorem B with non-trivial fundamental group for which the minimal entropy can be solved are the geometric ones studied by Paternain and Petean mentioned above. Furthermore if $M$ is simply connected and is a connected sum, then by another result of Paternain--Petean in \cite{PP} there exist only two such closed orientable $4$--manifolds which admit a metric of zero topological entropy, $\C {\rm P}^2 \# \C {\rm P}^2$ and $\C {\rm P}^2 \# \overline{\C {\rm P}}^2$. Modulo the case of $\mathbb{H}^2\times \mathbb{H}^2$ manifolds, these results provide a complete solution to the minimal entropy problem for geometrisable four--manifolds.

On the organisation of this paper: Section 2 contains definitions and the statement of theorem 1. Proposition 2 is shown in section 3. The proof of theorem 1 is a case by case analysis of the geometries in {\bf V} and it occupies sections 4 to 13. All the results are collected in section 14, where the proofs of theorems 1, A and B can be found.

{\bf Acknowledgements:} The results found here are part of my doctoral thesis. It is my pleasure to thank my supervisor Gabriel Paternain for his constant support, energetic help and kind encouragement. I also owe thanks to Dieter Kotschick, for his interest in the ideas presented here as well as for noticing a couple of mistakes and commenting on a previous version of this work. Likewise I would like to thank Ivan Smith and Jim Anderson for all the corrections and suggestions they made to my thesis, which have certainly improved the presentation here. I am also grateful to Alberto Verjovsky and Peter Scott, their experience and insight was a very positive influence. Finally I also wish to say thanks to Jimmy Petean and Dan Jane.

Financial support from the Consejo Nacional de Ciencia y Tecnolog\'ia, CONACyT M\'exico and the Secretar\'ia de Educaci\' on P\'ublica, SEP M\'exico during my time as a PhD student is gratefully acknowledged. Currently I am supported by the Deutsche Forschungsgemeinschaft, DFG Germany under the project "Asymptotic invariants of manifolds".

\section{Preliminaries}

 The simplicial volume $||M||$ of a closed orientable manifold is defined as the infimum of $\Sigma_{i}|r_{i}|$ where $r_{i}$ are the coefficients of a \emph{real} cycle representing the fundamental class of $M$.

For a closed connected smooth Riemannian manifold $(M,g)$, let ${\rm Vol}(M,g)$ be the volume of $g$ and let $K_{g}$ be its sectional curvature. We define the following minimal volumes as in \cite{Gro}:
\[ {\rm MinVol }(M):=\inf\limits_{g} \{ {\rm Vol}(M,g)\; : \; |K_{g}| \leq 1 \} \]
\[ {\rm Vol_{K} }(M):=\inf\limits_{g} \{ {\rm Vol}(M,g)\; : \; K_{g} \geq -1 \} \]

If ${\rm Vol_{K}}(M)$ is zero then the simplicial volume of $M$ is also zero, this follows from Bishop's comparison theorem \cite{PP}.

A $\mathcal{T}$--structure on a smooth closed manifold $M$ is a finite open cover $\{ U_{i} \}_{i=1}^k$ of $M$ with a non-trivial torus action on each $U_{i}$ such that the interesections of the open sets are invariant (under all corresponding torus actions) and the actions commute.
A $\mathcal{T}$--structure is called \emph{polarised} if the torus actions on each $U_{i}$ are locally free and in the intersections the dimensions of the orbits (of the corresponding torus action)  is constant. The structure is called \emph{pure} if the dimension of the orbits is constant.

By definition an $\Fs$--structure on a closed manifold $M$ is given by,
\begin{enumerate}
\item{ A finite open cover $\{ U_1, ..., U_{N} \} $;}
\item{ $\pi_{i}\co\widetilde{U_{i}}\rightarrow U_{i}$ a finite Galois covering with group of deck transformations $\Gamma_{i}$, $1\leq i \leq N$;}

\item{ A smooth torus action with finite kernel of the $k_{i}$-dimensional torus, \\ $\phi_{i}\co T^{k_{i}}\rightarrow {\rm{Diff}}(\widetilde{U_{i}})$, $1\leq i \leq N$;}

\item{ A homomorphism $\Psi_{i}\co \Gamma_{i}\rightarrow {\rm{Aut}}(T^{k_{i}})$ such that
\[ \gamma(\phi_{i}(t)(x))=\phi_{i}(\Psi_{i}(\gamma)(t))(\gamma x) \]
for all $\gamma \in \Gamma_{i}$, $t \in T^{k_{i}}$ and $x \in \widetilde{U_{i}}$; }

\item{ For any finite sub-collection $\{ U_{i_{1}}, ..., U_{i_{l}} \} $ such that  $U_{i_{1}\ldots i_{l}}:=U_{i_{1}}\cap \ldots \cap U_{i_{l}}\neq\emptyset$ the following compatibility condition holds: let $\widetilde{U}_{i_{1}\ldots i_{l}}$ be the set of points $(x_{i_{1}}, \ldots , x_{i_{l}})\in \widetilde{U}_{i_{1}}\times \ldots \times \widetilde{U}_{i_{l}}$ such that $\pi_{i_{1}}(x_{i_{1}})=\ldots = \pi_{i_{l}}(x_{i_{l}})$. The set $\widetilde{U}_{i_{1}\ldots i_{l}}$ covers $\pi_{i_{j}}^{-1}(U_{i_{1}\ldots i_{l}}) \subset \widetilde{U}_{i_{j}}$ for all $1\leq j \leq l$, then we require that $\phi_{i_{j}}$ leaves $\pi_{i_{j}}^{-1}(U_{i_{1}\ldots i_{l}})$ invariant and it lifts to an action on $\widetilde{U}_{i_{1}\ldots i_{l}}$ such that all lifted actions commute. }
\end{enumerate}

An $\Fs$--structure is said to be {\em pure} if all the orbits of all actions at a point, for every point have the same dimension. We will say an $\Fs$--structure is {\em polarised} if the smooth torus action $\phi_{i}$ above are fixed point free for every $U_{i}$.

The existence of a polarised $\Fs$--structure on a manifold $M$ implies the minimal volume ${\rm MinVol}(M)$ is zero, this is shown in \cite{CG}, the interested reader is invited to consult further examples there as well.

One of the main contributions of this paper is the proof of the following result.

\begin{Theorem}{  Let $M$ be a closed orientable geometric manifold modelled on a
    geometry in ${\bf{V}}$. Then $M$ admits a
    $\mathcal{F}$--structure.}
\end{Theorem}

   The case by case analysis which we will follow to prove this result was outlined in the introduction, a proof which collects all of these results can be found in the last section.

\section{Geometric manifolds of positive simplicial volume}

\begin{Proposition}\label{Hpos} {  If $M$ is a closed oriented geometric manifold modelled on a geometry in {\bf H}
     then  $\| M \| >0$.}
    \end{Proposition}
    \begin{proof}  First assume $M$ is closed, oriented and modelled in either
    $\mathbb{H}^{4}$ or  $\mathbb{H}_{\C}^{2}$.  Then $M$ admits a metric of negative sectional curvature and by the Thurston--Inoue--Yano theorem in \cite{Gro, IY} we have that  $\| M \| > 0$.  In fact this shows that any $M$ with \emph{finite volume} modelled on $\mathbb{H}^{4}$ and  $\mathbb{H}_{\C}^{2}$ has positive simplicial volume.

    Let $M$ be a closed manifold modelled on $\mathbb{H}^{2}\times\mathbb{H}^{2}$. We can use Gromov's Proportionality Principle from \cite{Gro} to see that $\| M \| > 0$.
     If the universal coverings of two closed Riemannian manifolds $M$ and $N$ are isometric then $\frac{\| M \|}{{\rm Vol}(M)}=\frac{\| N \|}{{\rm Vol}(N)}$.

     Consider the product of two hyperbolic surfaces $N = S_1 \times S_2$. The smooth manifold $N$ is modelled on $\mathbb{H}^{2}\times\mathbb{H}^{2}$. Because the simplicial volume of a product of closed manifolds is bounded from below by the product of their respective simplicial volumes we have  $\| N \|\geq C \| S_1 \| \| S_2 \| >0$ for some constant $C$ .

    Therefore  $\| M \| >0$ for any closed manifold $M$ modelled on $\mathbb{H}^{2}\times\mathbb{H}^{2}$. \end{proof}

\begin{rem} We expect complete manifolds of finite volume modelled on $\H^2 \times \H^2$ to have positive minimal entropy.
\end{rem}

\section{Circle Foliations}

\subsection{Geodesic foliations by circles and polarised $\Ts$--structures}

The following result explains how geodesic circle foliations can be used to define $\Ts$--structures on closed manifolds.

\begin{Proposition} Let $M$ be a closed manifold foliated by circles. Suppose $M$ admits
a metric such that the circles are geodesics. Then $M$ admits a polarised $\Ts$--structure.
\end{Proposition}

\begin{proof} By a theorem of A.W. Wadsley \cite{Wad}, the foliation by circles gives rise
to an orbifold bundle or Seifert fibration. This means that locally we have the following model
for the foliation near a fixed leaf $L$ (\cite[Theorem 4.3]{Ep1}).
There exists a finite group $G\subset {\rm O}(n)$ (where $\dim M=n+1$) and a homomorphism
$\psi \co \pi_{1}(L)=\Z\to {\rm O}(n)$. Let $\tilde{L}$ be the covering of $L$ corresponding to the
kernel of $\psi$, that $\tilde{L}$ is compact follows from \cite[Theorem 4.3]{Ep1}. Then $G$ acts on $\tilde{L}$ by deck tranformations
and we can consider the quotient $(\tilde{L}\times D^n)/G$, where $D^n$ is the unit ball
in $\re^n$.

Theorem 4.3 in \cite{Ep1} asserts the existence of $G$ and $\psi$ and a diffeomorphism
between $(\tilde{L}\times D^n)/G$ and a neighbourhood of $L$ preserving the leaves.

When $L$ is $S^1$, $\tilde{L}$ is also a circle and $G$ can only be a cyclic group $\Z_n$ and the obvious circle action
on $S^1\times D^n$ clearly descends to a circle action on $(S^1\times D^n)/G$ and thus locally
we always have a locally free circle action, whose orbits are precisely the leaves of the foliation.

If one can coherently orient all the leaves we would have a circle action on $M$, but if not, we still
have a $\mathcal T$-structure, since ``opposite'' actions still commute.
Let us make this a bit more precise.

The leaf space $B$ is an orbifold and $M$ is an orbifold bundle over $B$. As such, we may cover $B$ with compatible open sets
such that the transition maps of the bundle have values in ${\rm O}(2)$ (the fibres are circles and ${\rm Diff}(S^1)$ deformation retracts onto ${\rm O}(2)$).
But given $h\in {\rm O}(2)$ we obviously have $hlh^{-1}\in {\rm SO}(2)$ for any $l\in {\rm SO}(2)$.
Thus if we conjugate the obvious circle action of $S^1$ on itself by an element of ${\rm O}(2)$ we obtain
a new circle action commuting with the original one. Thus $M$ has a $\mathcal T$--structure. \end{proof}

Some of the four dimensional geometries are foliated by $\mathbb{R}$. This foliation descends to a circle foliation on their geometric manifolds, using this $S^1$-foliation we can define a $\mathcal{T}$--structure.

\begin{Theorem}\label{s1fol}
Every closed geometric manifold $M$ modelled on any of the geometries
\[ \mathbb{X}^{4}\in \{
\mathbb{S}^{3}\times \mathbb{E} , \mathbb{H}^{3}\times \mathbb{E} ,
\Sl\times \mathbb{E} , \mathbb{N}il^{3}\times \mathbb{E} ,
\mathbb{S}ol^{3}\times \mathbb{E}, \mathbb{N}il^{4}, \mathbb{S}ol^{4}_{1} \} \] admits a polarised $\mathcal{T}$--structure.

\end{Theorem}

 \begin{proof} In each of the geometries $\mathbb{S}^{3}\times \mathbb{E}, \mathbb{H}^{3}\times \mathbb{E},
\Sl\times \mathbb{E}, \mathbb{N}il^{3}\times \mathbb{E}$ or $
\mathbb{S}ol^{3}\times \mathbb{E}$ we have a trivial foliation given by the product with the Euclidean factor. In the case of $ \mathbb{N}il^{4} = \re^3 \ltimes_{\theta} \re$, with $\theta(t)=(t,t,t^2/2)$, it is given by the $\re$ factor on the right hand side of the semi-direct product.  For the remaining geometry of solvable Lie type
$$\mathbb{S}ol^{4}_{1}= \left\{
\left(
\begin{array}{ccc}
 1 & a  & c  \\
 0 &  \alpha &  b \\
 0 & 0  &  1
\end{array}
\right) : \alpha, a, b, c \in \re, \alpha > 0 \right\} $$
the $\re$ we are interested in is given by elements of the form
\[
\left(
\begin{array}{ccc}
 1 &  0 &  c \\
 0 &  1 &  0 \\
 0 & 0  & 1
\end{array}
\right) \quad \text{with} \quad c\in \re.
\]

This foliation on $\mathbb{X}^{4}$ descends to a foliation $\mathcal{F}$ on any quotient $M=\mathbb{X}^{4} / \Gamma$ under the action of a discrete group of isometries $\Gamma$ and the leaves of  $\mathcal{F}$ are all circles. The leaf space is  a geometric $3$--orbifold, with geometry $\mathbb{S}^{3}, \mathbb{H}^{3},
\Sl, Nil^{3},
\mathbb{S}ol^{3}, Nil^{3}$ and
$\mathbb{S}ol^{3}$ respectively.

By the proposition above $M$ admits a polarised $\mathcal{T}$--structure.
\end{proof}

A proof of existence of a polarised $\mathcal{T}$--structure on 4--manifolds which are foliated by 2--tori could be quite similar. We have an orbifold bundle
over a 2-orbifold with fibre $T^2$. If the transition maps may now be taken in ${\rm GL}(2,\Z)$, then the actions
defined in the local model when translated via the trivilizations to $M$ give rise
to commuting actions and we get again a $\mathcal T$--structure.
Actually ${\rm Diff}(T^2)$ deformation retracts onto the affine group of $T^2$ which is
$T^2\ltimes {\rm GL}(2,\Z)$ and this structure group will produce commuting actions.

\subsection{A general principle} Suppose $M\rightarrow B$ is a fibre bundle---or orbifold bundle---with fibre
$F$ and structure group $K\subset {\rm Diff}(F)$. Suppose the fibre $F$ admits a circle or torus action $\rho$
such that:
\begin{itemize}
\item it commutes with the action of the finite groups $G$ in the case of an orbifold;
\item if we conjugate $\rho$ by an element of $K$, the resulting action commutes with $\rho$.
\end{itemize}
Then $M$ admits a $\mathcal T$--structure.

We will verify this principle for the case of Seifert fibred smooth four manifolds in the next section.

\begin{rem} It will be a consequence of our construction of $\Fs$--structures on manifolds with decompositions into geometries of mixed type that some of the geometric manifolds we have just considered will also admit $\Fs$--structures (of course a $\Ts$--structure is already an $\Fs$--structure).
\end{rem}

\section{Seifert Fibred Geometries}

\subsection{Seifert Fibrations}

Let $S$ be a closed geometric manifold modelled on
$\mathbb{S}^{2} \times \mathbb{E}^{2}$,
$\mathbb{H}^{2}\times\mathbb{E}^{2}$,
$\Sl \times\mathbb{E}^{1}$, $\mathbb{N}il^{4}$,
$\mathbb{S}ol^{3}\times\mathbb{E}^{1}$ or
$\mathbb{S}^{3}\times\mathbb{E}^{1}$. It was shown by  M Ue in \cite{Ue, Ue2} that $S$ is a Seifert fibred space . We will briefly review the structure of Seifert fibrations in dimension four. We will review the description of these structures locally---in a neighbourhood of a point in $S$---this will allow us to furnish these manifolds with $\mathcal{T}$--structures. We refer to \cite{Ue} \& \cite{Ue2} throughout this section for details. The reader familiar with this data may want to skip to Lemma \ref{loc-ac} and then straight to Theorem \ref{seifert}, which are the main results of this section.

\begin{Definition}   A smooth oriented $4$--manifold $S$ is
Seifert fibred if it is the total space of an orbifold bundle $\pi \co S\rightarrow B$ with
general fibre a torus over a $2$--orbifold $B$.
\end{Definition}

Notice that the class of Seifert $4$--manifolds contains all the compact complex surfaces diffeomorphic to elliptic surfaces $X$ with $c_2(x)=0$ and also contains examples which do not admit any complex structure \cite{Wall}.

 Let $\pi \co S\rightarrow B$ be a Seifert fibration, with $S$ a geometric
manifold modelled in one of the geometries mentioned above and $B$
the orbifold base of the fibration. Denote by $T^2$ be the standard torus and $G\subset {\rm O (2)}$ a discrete subgroup, viewed as a group of Euclidean isometries. For any point $p\in B$ there exists a neighbourhood $U$ of $p$
such that $\pi^{-1}(U)$ is diffeomorphic to $(T^{2}\times D^{2})/G$ for some $G\subset {\rm O (2)}$. Here
$T^{2}$ is parametrised by two unit complex circles
$S^{1}\times S^{1}\subset \mathbb{C}^{2}$ , $D^{2}$
is the open unit complex disk $|z|\leq 1$ in $\mathbb{C}$ and $G$ is
the stabiliser at $p$, which acts freely on $T^{2}\times D^{2}$.

\subsection{Local description}

For $G$ non-trivial, there are three cases to consider, cyclic groups of rotations $ \mathbb{Z}_{m}$, reflection groups $\mathbb{Z}_{2}$ and dihedral groups $D_{2m}$.

\begin{enumerate}
    \item{ $G\cong \mathbb{Z}_{m}\cong \langle \rho \rangle$, where $\rho$ is a
    rotation of $\frac{2\pi}{m}$. This isotropy subgroup corresponds
    to cone points of cone angle $\frac{2\pi}{m}$. The action
    $\rho \co T^{2}\times D^{2}\rightarrow T^{2}\times D^{2}$ is given by
    $\rho(x,y,z)=(x-\frac{a}{m}, y - \frac{b}{m}, e^{\frac{2\pi i}{m}}z)$, here $x,y \in S^1$ and $z\in D^2$
    with $g.c.d.(m,a,b)=1$.
    The fibre over $p=0$ is called a multiple torus of type $(m,a,b)$.}

    \item{$G\cong \mathbb{Z}_{2}\cong \langle \ell \rangle$. Where $\ell$ is a
    reflection on the second factor of $T^{2}$ and on $D^{2}$.
    Now the action $\ell \co T^{2}\times D^{2}\rightarrow T^{2}\times D^{2}$ is given by
    $\ell(x,y,z)=(x+ \frac{1}{2}, -y, \bar{z})$. This is the isotropy
    subgroup corresponding to points on a reflector line or circle.
    In this case the fibre over $p$ is a Klein bottle $K$ and
    $\pi^{-1}(U)$ is a non-trivial $D^{2}$--bundle over $K$.}

    \item{$G\cong D_{2m} = \langle \ell, \rho : \ell^{2}=\rho^{m}=1,
    \ell\rho\ell^{-1}=\rho^{-1} \rangle $, for $m\in \mathbb{Z}$ . This is a dihedral group, the
    isotropy subgroup of corner reflector points of angle
    $\frac{\pi}{m}$, with the actions
    $\rho, \ell \co T^{2}\times D^{2}\rightarrow T^{2}\times D^{2}$ given by,
     \begin{eqnarray*}
    \rho(x,y,z) & = & (x , y - \frac{b}{m} , e^{\frac{2\pi i}{m}}z), \\
    \ell(x,y,z) & = & (x+ \frac{1}{2}, -y, \bar{z}).
    \end{eqnarray*}
      Informally, the fibre over $p$ is a Klein bottle whose fundamental domain is
    $\frac{1}{m}$--times that of the fibre of the reflector point near
    $p$. We call this fibre a multiple Klein bottle of type $(m,0,b)$. }
    \end{enumerate}

%
%

\subsection{ Local  $S^1$--actions} The descriptions above allow us to construct local circle actions on $S$.

\begin{Lemma}\label{loc-ac} Let $\pi \co S\rightarrow B$ be a Seifert fibration fo the $4$--manifold $S$.
For every point $p\in B$ there exists a neighbourhood $U$ of
$p$ diffeomorphic to  $D^{2}/G$ with $G \subset {\rm O}(2)$  such that $\pi^{-1}(U)\cong (T^{2}\times D^{2}) / G $ admits an $S^{1}$ action
which commutes with the action of $G$  on $T^{2}\times D^{2}$.
\end{Lemma}

\proof Take $U\subset B$ such that $\pi^{-1}(U)$ is diffeomorphic to
$(T^{2}\times D^{2})/G$, here we can define an $S^{1}$--action.
We will do this by first lifting the quotient by $G$ to $(T^{2}\times D^{2})$
and then showing that the $S^{1}$ action commutes in $(T^{2}\times D^{2})$
with the actions of all the different possible isotropy groups $G$. Hence
this $S^{1}$--action will be well defined in the quotient $(T^{2}\times
D^{2})/G \cong \pi^{-1}(U)\subset S$, this defines a local
$S^{1}$--action on $S$.

Define $\varphi \co S^{1}\times (T^{2}\times D^{2})\rightarrow (T^{2}\times
D^{2})$ by $\varphi(\theta, x,y,z)=(x+\theta, y, z)$. And now we
consider the various cases for $G$, the isotropy subgroup at $p\in U$.

 Let $G\cong \mathbb{Z}_{m}\cong  \langle \rho \rangle$, here $\pi^{-1}(U)\cong
(T^{2}\times D^{2})/ \langle \rho \rangle$. We claim $\varphi\circ
\rho=\rho\circ\varphi$:
\begin{eqnarray*}
\varphi\circ\rho & = & \varphi(x-\frac{a}{m}, y - \frac{b}{m} , e^{\frac{2\pi i}{m}}z
) \\
 & = & ((x-\frac{a}{m}) + \theta , y - \frac{b}{m} , e^{\frac{2\pi i}{m}}z) \\
 & = & ((x+ \theta) -\frac{a}{m}, y - \frac{b}{m} , e^{\frac{2\pi i}{m}}z) \\
 & = & \rho\circ\varphi
 \end{eqnarray*}

Consider $G\cong \mathbb{Z}_{2}\cong \langle\ell\rangle$, in this case $\pi^{-1}(U)\cong
(T^{2}\times D^{2})/\langle\ell\rangle$:
\begin{eqnarray*}
\varphi\circ\ell & = & \varphi(x+ \frac{1}{2}, -y, \bar{z}) \\
 & = & ((x+ \frac{1}{2}) + \theta , -y, \bar{z}) \\
 & = & ((x + \theta ) + \frac{1}{2}, -y, \bar{z}) \\
& =\ & \ell\circ\varphi
\end{eqnarray*}

When $G\cong D_{2m}$ it suffices to prove that
$\varphi$ commutes with both generators $\ell$ and $\rho$:
\begin{eqnarray*}
\varphi\circ\rho & = & \varphi(x , y - \frac{b}{m} , e^{\frac{2\pi i}{m}}z) \\
 & = & ( (x+ \theta) , y - \frac{b}{m} , e^{\frac{2\pi i}{m}}z)\\
 & = & \rho\circ\varphi. \\
\\
\varphi\circ\ell & = & \varphi(x+ \frac{1}{2}, -y, \bar{z})\\
& = &  ((x+ \frac{1}{2}) + \theta , -y, \bar{z}) \\
& = & ((x+\theta) + \frac{1}{2}, -y, \bar{z})=\ell\circ\varphi
\end{eqnarray*}\proved

So we now know that given an orbifold chart $U$ of $B$, we can construct an $S^1$--action on its preimage $\pi^{-1}(U)$  which is equivariant with respect to the action of $G$ on $T^2 \times D^2$.

\subsection{Description along the singular set}

The following picture along the reflector circles is taken from \cite{Ue} \& \cite{Ue2}.
Let $l$ an $h$ be the curves in $T^2$ represented by $\re / \Z \times \{ 0 \}$ and $ \{ 0 \} \times \re / \Z $ respectively. A choice of such a pair $(l,h)$ is called a framing for $T^2$.
The boundary of $B$ consists of a disjoint union of circles $C_i$, each of which we call a reflector circle. Let $N_i$ be an annulus bounded by $C_i$ and a curve $\gamma_i$ parallel to $C_i$.  In order to clarify the structure of $S$ near $C_i$ we now describe $\pi^{-1}(N_i)$. Say the corner reflectors $p_1, ... , p_s$ on $C_i$ are of type  $(m_1, 0, b_1), ... , (m_s, 0 , b_s)$, with respect to the framing $(l_i, h_i)$ of the general fibre over some base point of $N_i$.

Understanding the fibres over these corner reflectors is simplified if we consider the double cover $DB$ of $B$, with the projection $p \co DB\rightarrow B$ obtained by identifying 2 copies of $B$ along the reflector circles. Let $DN_i$ be the suborbifold of $DB$ covering $N_i$ and $D\pi \co DS\rightarrow  DB$ be the fibration induced from $\pi \co S \rightarrow B$. Then $S$ is the quotient of $DS$ by a free involution $\iota$ which is the lift of the reflection $l$ which switches both copies of  $DB$. The action of $\iota$ on the reflection point near the base point is identical to that of $l$ in case (3) above.
In the presentation of $\pi_1(S)$ in \cite{Ue} the map $\iota$ satisfies
\[ \iota^2=l \quad \text{and} \quad \iota h \iota^{-1}=h^{-1}. \quad \quad \quad\quad\quad\quad( \star ) \]
The corner reflector point $p_j$ is covered by  a cone point $q_j \in DB$ and the fibre over $q_j$ is a multiple torus of type $(m_j, 0, b_j)$.

\begin{figure}[htbp]
\begin{center}
\includegraphics[width=70ex]{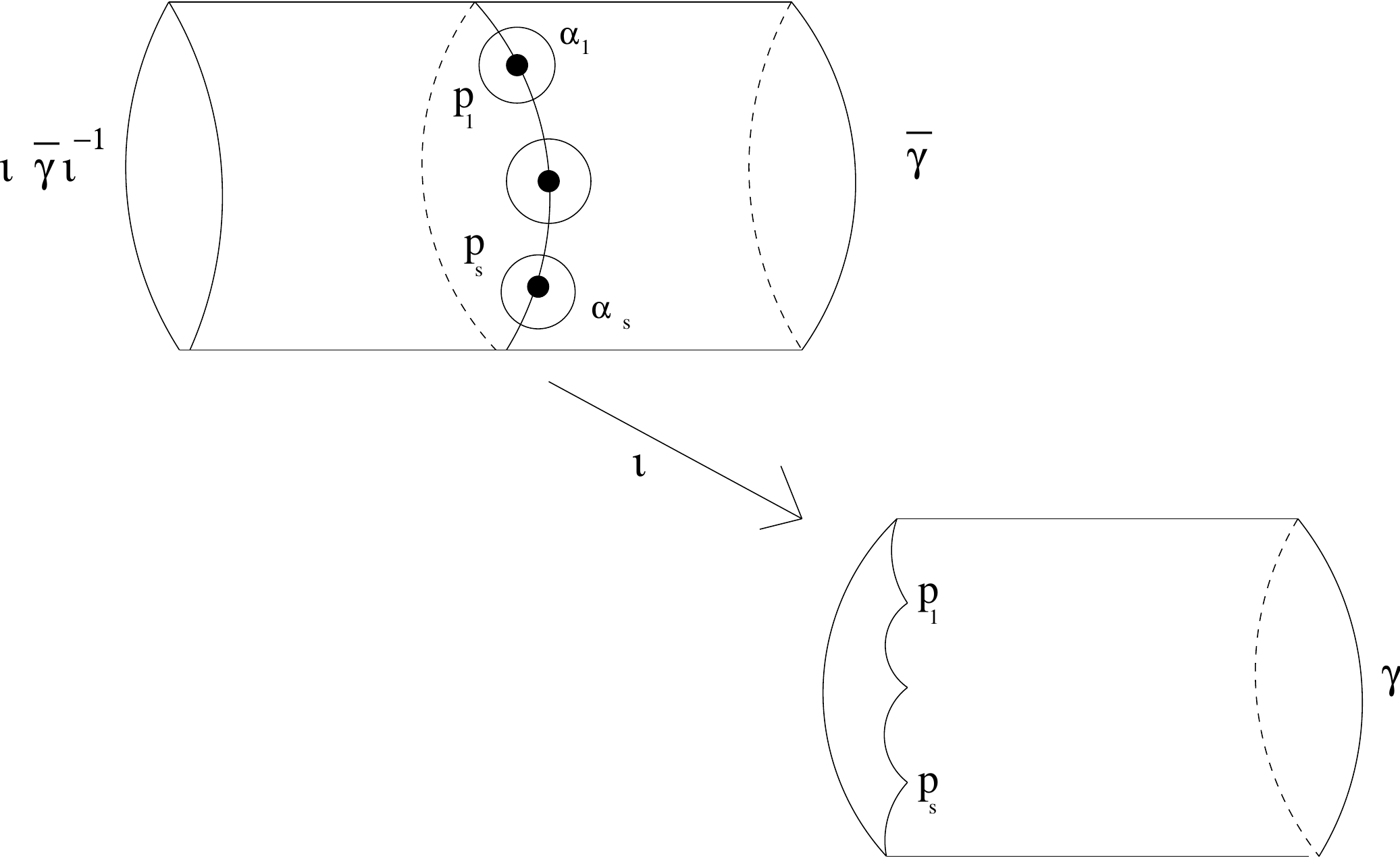}

\caption{Local picture along a reflector circle.}

\end{center}
\end{figure}

Take the oriented meridional circle $\alpha_{j}$ centered at $p_{j}$ as in figure 2, then the lifts $\tilde{\alpha}_1 , \dots, \tilde{\alpha}_{s}$ of the curves $\alpha_1, \ldots, \alpha_{s}$ can be taken to satisfy in $\pi_1(S)$ the following relations:
$$\tilde{\alpha}_{i}^{m_{i}} h^{b_{i}}=1 \quad (i=1, \ldots , s), \quad \iota \tilde{\alpha}_{s}\iota^{-1}= \tilde{\alpha}_{s}^{-1},$$
$$ \iota \tilde{\alpha}_{s-1}\iota^{-1}= \tilde{\alpha}_{s}^{-1} \tilde{\alpha}_{s-1}^{-1}  \tilde{\alpha}_{s}, \quad\ldots  , \quad \iota \tilde{\alpha}_{1}\iota^{-1}= \tilde{\alpha}_{s}^{-1}\tilde{\alpha}_{s-1}^{-1}\ldots \tilde{\alpha}_{1}^{-1}  \tilde{\alpha}_{2} \ldots  \tilde{\alpha}_{s} .$$

We can describe the monodromy of the fibration along a reflector circle. Let $V$ denote the union of small disk neighbourhoods around each corner reflector points $p_{j}$ and take $\bar{\gamma}$ and $\iota \bar{\gamma} \iota^{-1}$ as in figure 2. Then the curve represented by $ \bar{\gamma}^{-1} \alpha_1 \alpha_2 \ldots \alpha_s \iota \bar{\gamma} \iota^{-1}$ is null-homologous in $DN - V$. Hence the monodromy matrix $A$ along $\bar{\gamma}$ with respect to the framing $(l,h)$ must satisfy $JAJ^{-1}=A$ where $J=\begin{pmatrix}
     1 & 0   \\
     0 & -1 
\end{pmatrix} $ is a monodromy matrix for $\iota$. Then we must have that $A=\pm I$ , where $I$ is the identity matrix.

Another piece of information that we need in order to describe the Seifert fibration is the obstruction to extending the fibration over a neighbourhood of each reflector circle $C_{i}$, this is called the \emph{euler class} of $C_{i}$. Consider a lift $\tilde{\gamma}$ of $\gamma$, which also determines a lift $\iota \tilde{\gamma} \iota^{-1}$. Then we have the relation $ \tilde{\gamma}^{-1}\tilde{\alpha}_1 \ldots \tilde{\alpha}_{s} \iota \tilde{\gamma} \iota^{-1}=l^{a}h^{b}$ in $\pi_1(S)$. The euler class $(a,b)$ of $C_{i}$ is the obstruction to extending $\gamma \cup \iota \gamma \iota^{-1} \cup \alpha_1 \cup \ldots \cup \alpha_{s}$ to the cross section on $\pi^{-1}(DN - V)$. We have that $a=-1 $ if the monodromy $A$ around $\gamma$ is $-I$ and $a=0$ if $A=I$. The value of $b$ depends on the choices of the lifts $\tilde{\gamma} ,\tilde{\alpha_{i}}$ of $\gamma , \alpha$.

\subsection{Global description of a Seifert Fibration}

We can now give a global description of our Seifert fibred manifold $S$.

Let $N_{i}$ be a tubular neighbourhood of each reflector circle $C_{i}$, with boundaries $C_{i}$ and $\gamma_{i}$ as in the figures above. Fix a base point near $C_{i}$ and a framing $(l,h)$ of the general fibre satisfying $(\star )$. Denote by $|B|$ the topological space underlying the orbifold $B$.  Let $p_{i}$ be a cone point and $D_{i}$ a disk neighbourhood of $p_{i}$. If we fix the lift $\tilde{\gamma}_{i}$  of $\gamma_{i}$, then the fibration over the complement $B_{0}=B-\cup_{i}N_{i}$ is described by the following information;

(1) The monodromy matrices $A_{i}, B_{i}\in {\rm SL}(2, \Z)$ along the set of standard generators $s_{i}, t_{i}$ (for $i=1, \dots , g)$ of $\pi_1 |B_{0}|$ if $|B_{0}|$ is orientable.

(1') The monodromy matrices $A'_{i}\in {\rm GL}(2, \Z)$ with $\det A_{i}'=-1$ along the set of standard generators $v_{i}$ (for $i=1, \dots , g)$ of $\pi_1 |B_{0}|$ if $|B_{0}|$ is non-orientable.

(2) The type $(m_{i}, a_{i}, b_{i})$ of the multiple torus over the cone point $p_{i}$, (for $i=1, \dots , t)$.

(3) The obstruction $(a',b')$ to extending the fibration over $(\cup \gamma_{i})\cup (q_{i}')$ to the cross section in $\pi^{-1}(B_{0}\cup D_{i})$, where $q_{i}'$ is the lift of a meridional circle centred at $p_{i}$. This is called an euler class.
%
%
%

The fibration over $N_{i}$ is described as before with respect to the framing $(l_{i}, h_{i})$ of the general fibre on $N_{i}$ and the lift $\iota_{i}$ of the reflection along $C_{i}$ (where $(l_1, h_1)=(l,h)$) satisfying $\iota_{i}^2=l_{i}$ and $\iota_{i} h_{i} \iota^{-1}_{i}=h^{-1}_{i}$.

Then $\pi_1(N_{i})$ is attached to $\pi^{-1}(B_{0})$ so that $(l_{i}, h_{i})= (l,h) P$ for some $P\in {\rm SL}(2, \Z )$ with $P_1 = Id$. This implies that if we take the lift $\tilde{\delta}_{i}$ of the curve $\delta_{i} $ as in figure 3, then the monodromy along $\tilde{\delta}_{i}$ is  $B_{i}=P_{i}JP_{i}^{-1}J$ with respect to $(l,h)$. It is possible to take $\tilde{\delta}_{i}$ so that $\iota_{i}=\tilde{\delta}_{i}\iota$, because $\iota_{i}\tilde{\delta}_{i}\iota^{-1}=\tilde{\delta}_{i}^{-1}$ and so $\tilde{\delta}\iota_{i}\tilde{\delta}_{i}\iota^{-1}= l^{s+1}h^{t}$ for some $s,t \in \Z$ (recall that $\iota ^{2}=l$).

As a final step, we describe the relations between the monodromies. Let $I_{i}=\pm Id$ be the monodromy along $C_{i}$, and $A_{i}, B_{i}$ ( or $A'_{i}$ in case $B_0$ is not orientable ) be the monodromies along the standard curves on $B_0$ as before. Then $\prod [  A_{i}, B_{i} ] \prod  I_{i}=Id$ ( or $\prod  A_{i}^2\prod  I_{i}=Id$). The Seifert fibration of $S$ is determined by the information above, using this description we can now show the existence of $\Ts$--structures on Seifert fibred four manifolds.

\subsection{Seifert Fibrations and polarised $\Ts$--structures}

\begin{Theorem}\label{seifert}
 Every smooth closed and oriented Seifert fibred four manifold $S$ admits a polarised  $\mathcal{T}$--structure.
\end{Theorem}
\begin{proof} Let $\pi \co S\rightarrow B$ be a Seifert fibred smooth $4$--manifold,
over the orbifold $B$. So $S$ is the total space of an orbifold bundle with general
fibre a torus, over the $2$--orbifold $B$. Let $N_{i}$ be open
annular neighbourhoods of the circle reflectors $C_{i}$ of $B$. Take
$B_{0}=B - \bigcup_{i=1}^{r}N_{i}$.

Let $\mathcal{U}$ be an open covering for $B_{0}$ such that for $U $ in
$\mathcal{U}$ we have that  $\pi^{-1}(U)\cong (T^{2}\times D^{2})/G$ as in Lemma \ref{loc-ac} above. So $G$ is either
trivial or isomorphic to $\mathbb{Z}_{p}$. As $B_{0}$ is compact we
may choose a finite subcovering $\{ U_{i} \} $ of $\mathcal{U}$.

Notice that for $G$ trivial we have that
$(T^{2}\times D^{2})/G =
T^{2}\times D^{2}$
and for $G\cong \mathbb{Z}_{p}$, we have that $(T^{2}\times
D^{2})/G$
is diffeomorphic to $T^{2}\times D^{2}$. 

Denote $\pi^{-1}(B_0)$ by $S'$, the restriction $\pi |_{S'}\co S'\rightarrow B_{0}$ is an orbifold bundle with fibre $T^2$. The singular points of $B_0$ are all cone points.

For each $p\in U_{i}\cap U_{j}$ call the local trivialisations
$\Phi_{i}\co \pi^{-1}(U_{i})\rightarrow T^{2}\times D^{2}$ and $\Phi_{j}\co\pi^{-1}(U_{j})\rightarrow T^{2}\times D^{2}$. These give rise to the transition functions
$\Phi_{j}\circ \Phi_{i}^{-1}(x,y,z)=(\Psi_{ij}(z),z).$

Let $(l,h)$ be a framing for general fibre $T^{2}\cong
\mathbb{R}^{2}/\mathbb{Z}^{2}$, as explained above. Because the diffeomorphism group ${\rm Diff}(T^2)$ of $T^2$ retracts to $T^2 \ltimes {\rm GL}(2, \Z)$, the transition functions $\Psi_{ij}(z)$ can be
regarded as elements of  $ {\rm GL}(2,\mathbb{Z})$. So the
structural group of the orbifold bundle reduces to this linear one.

Showing that the actions $\varphi_{i}$ and $\varphi_{j}$  commute in the
intersections $\pi^{-1}(U_{i})\cap \pi^{-1}(U_{j})$ is an easy exercise in linear algebra. Therefore they define  a $\mathcal{T}$--structure on $\pi^{-1}(B_{0})$.

Now we exhibit a $\mathcal{T}$--structure on the neighbourhoods $N_{i}$ of the
reflector circles $C_{i}$. Consider $\overline{N}_{i}$, the closure
of $N_{i}$, as $N_{i}$ is covered by  open subsets in which we
defined the circle actions we extend these to cover $\overline{N}_{i}$ and take a finite subcovering $\{ V_{k} \}$ by. Let $\gamma_{i}$ denote the boundary of
$N_{i}$ which is not $C_{i}$. We claim the corresponding actions in the
sets $U_{j}$ in $B_{0}$ and $V_{k}$ in $N_{i}$ commute in the
intersection of these sets.

This follows from the fact that the fibration on $N_{i}$ is described
with respect to the framing $(l_{i},h_{i})$ of the general fibre on
$N_{i}$, as we will now see. Recall that $\pi^{-1}(N_{i})$ is attached to $\pi^{-1}(B_{0})$ so
that $(l_{i},h_{i})=P_{i}(l,h)$ for some $P_{i}\in  {\rm SL}(2, \mathbb{Z} )$, following the global description of a Seifert bundle above. Once more the matrices involved here behave well with respect to the actions we defined on $\pi^{-1}(N_{i})$ and $\pi^{-1}(B_{0})$ , meaning that the actions commute in the
intersection of these sets.

Therefore we have a $\mathcal{T}$--structure on the Seifert fibred
$4$--manifold S.
\end{proof}
\begin{Corollary} Every closed geometric four manifold $M$ modelled on one of the geometries
$\mathbb{S}^{2} \times \mathbb{E}^{2}$,
$\mathbb{H}^{2}\times\mathbb{E}^{2}$,
$\Sl\times\mathbb{E}^{1}$,  $\mathbb{N}il^{4}$,
$\mathbb{S}ol^{3}\times\mathbb{E}^{1}$ or
$\mathbb{S}^{3}\times\mathbb{E}^{1}$ admits a polarised $\mathcal{T}$--structure.
\end{Corollary}

\section{Flat four manifolds}\label{flat}

There are 27 orientable and 48 non-orientable compact flat 4-manifolds, out of which only 8 admit a complex structure \cite[p.199]{BHPV} and only 3 are not Seifert fibred \cite[p.146]{Hil}.

\begin{Lemma}
Every closed flat four manifold admits a polarised $\Fs$--structure.
\end{Lemma}

\begin{proof} This follows from a theorem of Bieberbach which states that every flat manifold is finitely covered by a torus of the same dimension. This covering constitutes a polarised $\Fs$-structure, as the torus acts freely on itself.
\end{proof}

\section{Solvable Lie Geometries.} We have already dealt with the geometries $\mathbb{S}ol^{4}_{1}$ and $\mathbb{S}ol^3\times \mathbb{E} =  \mathbb{S}ol^{4}_{n,n}$ in Theorem \ref{s1fol} where we gave manifolds modelled on them a locally free $S^1$-action. We will now focus on the remaining cases.

Recall that if we glue two manifolds along components of their boundary using  isotopic diffeomorphisms, the resulting manifolds are diffeomorphic, see for example \cite{Hir}.

\begin{Theorem}\label{maptorus} If $M$ is an orientable geometric manifold modelled on
    $\mathbb{S}ol^{4}_{0}$ or $\mathbb{S}ol^{4}_{m,n}$ when  $(m\neq n)$,
    then $M$ admits a polarised $\Ts$--structure.
\end{Theorem}

\begin{proof} Firstly we recall that if $M$ admits one of these
two geometries then Hillman has shown $M$ is homeomorphic to the mapping torus $M_{f}$ of a homeomorphism of
$T^{3}=\mathbb{R}^{3} / \mathbb{Z}^{3}$ \cite{Hil}. Moreover, Cobb \cite[p.176]{Cb} extended this to show that if $M$ is orientable then it is in fact diffeomorphic to such a mapping torus. Hence $M$ is diffeomorphic to  $M_{f} = (T^{3}\times I)/ \mathbb{s}$ , where  $(x,0)\mathbb{s} (f(x),1)$ for
some diffeomorphism $f$ of $T^{3}$.
Furthermore we also know
that its fundamental group is $\mathbb{Z}^3\rtimes_{A}\mathbb{Z}$ for
some $A\in  {\rm SL}(3,\mathbb{Z})$, as we are assuming $M$ is orientable.
We know by corollary \ref{3diff}  that any diffeomorphism of $T^3$ is isotopic to an affine transformation, see also \cite{Iv}. In this case $f$ is isotopic to the transformation induced by $A$ on $T^3$. Denote by $M_{A}$ the mapping torus of $T^3$ under the transformation induced by $A$.
That  $M_{f}$ is diffeomorphic to $M_{A}$ follows from the fact that mapping tori of isotopic diffeomorphisms are diffeomorphic.
Let $\varphi_t$ denote the action of $T^{3}$ on  $T^{3}\times \{ t \} $ by translations, and $\overline{\varphi_t}$ the lift to $\mathbb{R}^{3}$.
In order to define a $\mathcal{T}$--structure on $M_{A}$  using the actions  $\varphi_t$, we must now verify that $\varphi_{0}$ commutes with $\varphi_{1}$ when conjugated by the transformation induced by $A$ on $T^3$. A simple calculation (which is a particular case of Lemma \ref{nilaff} below) shows
$$A^{-1}\circ \overline{ \varphi_{1}} \circ A \circ \overline{ \varphi_{0}}= \overline{\varphi_{0}}\circ A^{-1}\circ \overline{\varphi_{1}} \circ A. $$
 Therefore $\varphi_{0}$ commutes $\varphi_{1}$ on  $M_{A}$. Notice that the dimension of every orbit is $3$ and that the action  $\varphi_t$  is locally free. Hence we have endowed $M_{A}$, and therefore $M$,  with a polarised $\mathcal{T}$--structure.  \end{proof}

\begin{rem}
The same argument constructs a polarised $\Ts$--structure on mapping tori of diffeomorphisms of $T^{n}$ which are isotopic to an affine transformation.
\end{rem}

\section{Sphere foliations} We will now see that when a compact closed manifold $M$ is modelled on a geometry of type $\mathbb{S}^2\times \mathbb{X}^2$, where $\mathbb{X}^2$ is a $2$--dimensional geometry, $M$ admits a $\mathcal{T}$--structure.

\begin{Theorem} A smooth orientable closed $4$--manifold which is foliated by $S^2$ or
  $\mathbb{R} {\rm P}^2$ admits a  $\mathcal{T}$--structure.
  \end{Theorem}

\begin{proof} Suppose $M$ is a smooth orientable closed $4$--manifold which admits a codimension 2 foliation
with leaves $S^2$ or $\mathbb{R} {\rm P}^2$. If all the leaves are
 homeomorphic then the projection to the leaf space is a submersion and $M$ is the total space of an $S^2$ or an
 $\mathbb{R} {\rm P}^2$ bundle over a surface. Then $M$ is known to admit effective $S^1$--actions (see Melvin \cite{Mel} and Melvin--Parker \cite{MP}). Assume that the leaves are not all homeomorphic.
 Having such a foliation is equivalent to having an $S^2$ orbifold
 bundle over a 2-orbifold (see  Ehresman \cite{Eh}, Epstein \cite{Ep1}, Eells--Verjovsky \cite{EV} and Molino \cite{Mol}), as such a foliation is Riemannian.

 Denote by $F$ the orbit space of the foliation and
 $\pi \co M\rightarrow F$ the orbifold bundle. Ehresman's structure theorem \cite{Eh, Ep1} implies its singularities may only
 be isolated points and provides the following description; for any point $p\in F$ and a small neighborhood $U$
 of $p$, $\pi^{-1}(U)$ is diffeomorphic to $(S^2\times D^2)/G$. Here
 $G$ is a discrete subgroup of $ {\rm O}(2)$ which acts freely on $S^2\times
 D^2$. Because $\mathbb{Z}_2$ is the only such group that acts freely on
 $S^2$ the only  possible
 singularites for the orbifold bundle correspond to projective planes
 $\mathbb{R} {\rm P}^2$ over the set of singular points $p_{i}\in F$.

Consider an open neighborhood $V_{i}$ of $p_{i}$, let $V=\cup
V_{i}$. Then the restriction of $\pi$ to $E=F-V$ is a fibre bundle with
total space $N\subset M$ and fibres $S^2$. Melvin and Parker  have shown
that $N$ admits an $S^1$ action given by rotations in the
fibres \cite{Mel} \& \cite{MP}. Moreover, they show that the structure group of $N\rightarrow E$ is contained in $ {\rm O}(2)$. Since ${\rm Diff}(S^2)$ retracts to $ {\rm O}(3)$ and preserves fibres, the transition maps are either isotopic to the identity or the antipodal map.

We
will now show that this is also the case for $M\rightarrow F$, this will
allow us to extend this action to a $\mathcal{T}$--structure on $M$.
Let $r_{\alpha}$ denote the rotation of $S^2$ with respect to the axis $\alpha$ and $a$ the antipodal map. An easy exercise in linear algebra shows these two transformations commute, that is $r_{\alpha}\circ a = a \circ r_{\alpha}$. Therefore the following diagram commutes.
$$\begin{array}{ccc}
S^2 & \stackrel{r_{\alpha}}\longrightarrow & S^2  \\
\downarrow  & & \downarrow \\
\mathbb{R} {\rm P}^2 & \stackrel{r_{[\alpha]}}\longrightarrow &  \mathbb{R} {\rm P}^2
\end{array}$$

Where  $r_{[\alpha]}$
  denotes the rotation of $\mathbb{R} {\rm P}^2$ with fixed point the class
  of $\alpha$.

For a neighborhood $V_{i}$ of a singular point $p_{i}$, we can
lift the preimage $\pi^{-1}(V_{i})$ to $S^2\times D^2$. The
action of $r_{\alpha}$ on $S^2\times D^2$ commutes with the quotient
of $\mathbb{Z}_2$, thus defining an $S^1$ action on $(S^2\times
D^2)/\mathbb{Z}_2$ which is diffeomorphic to $\pi^{-1}(V_{i})$.

The holonomy around $\partial V_{i}$ is $\mathbb{Z}_2$, so that the maps that attach $\pi^{-1}(V_{i})$ to
$N$ in order to obtain $M$ are either isotopic to the identity or to
the antipodal map. In the case of the identity there is nothing to
prove. If the attaching map is isotopic to the antipodal map it suffices to note that the rotations on $S^2$ which are defined on $N$ and $\pi^{-1}(V_{i})$ both commute with the antipodal map $a$. Therefore they define a $\mathcal{T}-$structure on $M$.
\end{proof}

Conveniently enough, Hillman has shown that if a manifold $M$ admits a geometric decomposition into pieces modelled on geometries of the type $\mathbb{S}^2\times \mathbb{X}^2$ then $M$ is foliated by $S^2$ or
  $\mathbb{R} {\rm P}^2$ \cite{Hil}. We use this description to see that we have also proved the following two results.

\begin{Corollary} Any smooth orientable 4-manifold $M$ with a geometric decomposition into pieces of the type $\mathbb{S}^2\times \mathbb{X}^2$ admits a $\mathcal{T}$-structure and hence has zero
    minimal entropy.
    \end{Corollary}

    In particular if $M$ is an orientable geometric four manifold modelled in a geometry of type $\mathbb{S}^2\times \mathbb{X}^2$ then $ {\rm h} (M)=0$.

\begin{Theorem}\label{s2foliations} A closed manifold $M$ modelled on a geometry of type $\mathbb{S}^2\times \mathbb{X}^2$, where $\mathbb{X}^2$ is a 2-dimensional geometry,  admits a $\mathcal{T}$-structure.
\end{Theorem}

    In these cases it is the best we can hope for, in general such a manifold $M$ might not admit a \emph{polarised} $\mathcal{T}$--structure because $M$ could have positive Euler characteristic $\chi (M)>0$ and therefore its minimal volume could not vanish.

\section{Geometric Decompositions}



\begin{Definition}{
We say that an $n$-manifold $M$ admits a geometric decomposition if it has a finite collection of disjoint 2-sided hypersurfaces $S$ such that each component of  $M -\cup S$ is geometric of finite volume.  }
\end{Definition}

In other words, each component of $M -\cup S$ is homeomorphic to $X / \Gamma $, for some geometry $\mathbb{X}$ and a lattice $\Gamma$. We shall call the hypersurfaces $S$ {\emph{cusps}} and the components of $M-\cup S$ {\emph{pieces}} of $M$. The decomposition is {\emph{proper}} if the set of cusps is nonempty.

\subsection{Dimension four}

In dimension four Hillman \cite[p.138]{Hil} brought together various
results and organised them to show that if a closed $4$--manifold $M$ admits a geometric decomposition
    then either,
{\begin{enumerate}
\item{$M$ is geometric;}
\item{$M$ has a codimension 2 foliation with leaves $S^2$ or $\re {\rm P}^2$;}

\item{the pieces of $M$ have geometry $ \mathbb{H}^{4}$,
  $\mathbb{H}^{3}\times \mathbb{E}$, $ \mathbb{H}^{2}\times
  \mathbb{E}^{2}$ or $ \Sl \times \mathbb{E}$;}

\item {the pieces of $M$ have geometry $\mathbb{H}^{2}_{\mathbb{C}}$ or $ \mathbb{F}^{4}$;}
  \item{the pieces of $M$ all have geometry
  $\mathbb{H}^{2}\times \mathbb{H}^{2}$.}
\end{enumerate}}

In the first 3 cases $\chi(M)\geq 0$, in the last 2 cases $M$ is aspherical.

This follows from inspecting the various possible types of cusps that
appear in a geometric decomposition.  We synthesise this
information in the following table where we collect the information about the geometric structures on the cusps of geometric four manifolds, as well as the references where these results can be found.

\xymatrix{ {\rm Geometry} & {\rm Cusps} & {\rm Reference} \\
\mathbb{H}^{n} & {\rm flat} & {\rm Eberlein \cite{Eb}} \\
\mathbb{H}^{3}\times \mathbb{E} \;  \mathbb{H}^{2}\times
  \mathbb{E}^{2}\; {\rm and}\;  \Sl \times \mathbb{E} & {\rm flat} & {\rm Hillman  \cite{Hil}} \\
\mathbb{S}^{2}\times \mathbb{H}^{2} &  \mathbb{S}^{2}\times \mathbb{E}{\rm -manifolds} & {\rm Hillman \cite{Hil}} \\
\mathbb{F}^{4} &  \mathbb{N}il^{3}{\rm -manifolds} & {\rm Hillman \cite{Hil}} \\
\mathbb{H}^{2}_{\mathbb{C}} & \mathbb{N}il^{3}{\rm -manifolds} & {\rm Goldman \cite{Gld}}}


The cusps of irreducible $ \mathbb{H}^{2}\times
  \mathbb{H}^{2}$ manifolds are modelled on $Sol^{3}$, they
  are graph manifolds whose fundamental group contains a non-abelian
  subgroup otherwise \cite{Shi}.

These are the only geometries we need to consider, because if a
geometry is of solvable or compact type every lattice has compact quotient \cite{Ra}.

\subsection{Geometrisable 4-manifolds and positive simplicial volume}

A manifold is called {\emph{geometrisable}} if it is diffeomorphic to a connected sum of manifolds which admit geometric decompositions. Given a manifold $N$ with a geometric decomposition, if the fundamental group of the hypersurfaces of the decomposition injects into $\pi_1(N)$ we say that the geometric decomposition is $\pi_1$--injective.

In dimension 4, Hillman observed that, except for reducible pieces modelled on the geometry $\H^2 \times \H^2$, the inclusion of a cusp into the closure of a piece induces a monomorphism on the fundamental group. That is to say (again, except for reducible $\H^2 \times \H^2$-pieces ) every geometric decomposition in dimension four is $\pi_1$--injective \cite[p.139]{Hil}.

Therefore we can show the following holds:

\begin{Proposition}\label{Hdecpos} { Let $M$ be a geometrisable smooth four manifold. If a piece of the decomposition of $M$ is modelled on the real hyperbolic four-space $\H ^4$ or on the complex hyperbolic plane $\H^2_{\C } $ then $\| M \| >0$.}
\end{Proposition}
\begin{proof} Let $N$ denote a piece of $M$ modelled on $\H ^4$ or $\H^2_{\C } $. The manifold $N$ has finite volume and negative curvature, which implies $\| N \| >0$.

The cusps of manifolds modelled in the geometries $\H ^4$ and $\H^2_{\C } $ are either flat or $Nil^3$-manifolds, respectively.

This implies that we can cut the cusps $S$ of $N$ off from $M$ (by Gromov's Cutting--Off theorem \cite[p.58]{Gro}), because the fundamental group of any cusp of $N$ is amenable. This means the simplicial volume of $M$  is not affected when we take the cusps of $N$ off. In other words, if $N^0 := M -(N \cup S)$ then the cuttuing off theorem implies $\| M \| = \| N^0 \| + \| N \|$.

Therefore $\| M \| \geq \| N\| >0$.
\end{proof}

We expect the same to hold true for geometric decompositions with pieces modelled on the product of two hyperbolic planes $\H^2 \times \H^2$.

\begin{Conjecture}
If a smooth orientable four manifold $M$ admits a proper geometric decomposition into pieces modelled on $\H^2 \times \H^2$ then ${\rm h}( M)  \neq 0$.
\end{Conjecture}

\section{Mixed euclidean cases $\H^3\times \E$, $\H^2\times \E^2$ and $\Sl\times\E$}

In this section we will deal with the manifolds in case $(3)$ of Hillman's Theorem which do not have pieces modelled on $\H^4$.

\subsection{Generalities on the isometry group of $\mathbb X$}
\begin{Definition} A Riemannian manifold $M$ is reducible if $M$ is isometric to the Riemannian product $M_{1}\times M_{2}$ of two manifolds, $M_{1}$ and $M_2$ of positive dimension. If $M$ is not reducible, then it is said to be {\bf irreducible}.
\end{Definition}

In general if we have a simply connected Riemannian product $N\times M$, where $M$ is Euclidean
and $N$ is irreducible (a de Rham decomposition in the notation of \cite{Eb}), then
${\is}(N\times M)={\is}(N)\times{\is}(M)$ (see \cite[p. 240]{KN}). Thus
$\is(\H^3\times\E)=\is(\H^3)\times\is(\E)$, $\is(\H^2\times\E^2)=\is(\H^2)\times\is(\E^2)$
and $\is(\Sl\times\E)=\is(\Sl)\times\is(\E)$.

The identity components of these groups are:
\begin{align*}
&\is_{0}(\H^3\times \E)=PSL(2,\C)\times \re,\\
&\is_{0}(\H^2\times \E^2)=PSL(2,\re)\times \is^{+}(\E^2),\\
&\is_{0}(\Sl\times \E)=\is_{0}(\Sl)\times \re.
\end{align*}
The group $\is(\Sl)$ has only two connected components and no orientation reversing
isometry \cite{Scot}.

\subsection{Lattices}\label{lat} Let $\Gamma\subset \is(\X)$ be a discrete subgroup
which acts freely on $\X$ such that $M:=\X/\Gamma$ is a complete orientable manifold
with finite volume.

By a theorem of Wang (cf. \cite[8.27]{Ra}),
the lattice $\Gamma$ meets the radical $R$ of the connected Lie group
$\is_{0}(\X)$ in a lattice. The radicals are Euclidean and may be described as follows.
For $\H^3\times \E$, the radical
is the copy of $\re$ given by the translations on the $\E$ factor.
For $\H^2\times \E^2$, the radical is a copy of $\re^2$ given by the translations on the $\E^2$ factor.
For $\Sl\times \E$ it is also $\re^2$, with one copy of $\re$ coming from translations
on the $\E$ factor and the other coming from the center of $\is_{0}(\Sl)$.
Thus $\Gamma\cap R$ is isomorphic to $\Z$ or $\Z^2$.

\subsection{$\Fs$--structures on flat manifolds}\label{flatend}
The isometry group of $\E^n$ is the semidirect product of $\re^n$ and $ {\rm O}(n)$.
Let $\rho\co  {\rm O}(n)\to \mbox{\rm Aut}(\re^n)$ be the map $\rho(B)(x)=Bx$.
Let $\Gamma\subset \is(\E^n)$ be a cocompact lattice and $M:=\E^n/\Gamma$ a compact flat manifold.
Let $p\co \Gamma\to  {\rm O}(n)$ be the homomorphism $p(t,\alpha)=\alpha$, where $(t,\alpha)\in \re^n\times  {\rm O}(n)$.
The Bieberbach theorem guarantees that $\Gamma$ meets the translations in a lattice
(ie the kernel of $p$ is isomorphic to $\Z^n$) and $p(\Gamma)$ is a finite group $G$.
Then $M$ is finitely covered by the torus $\re^n/\mbox{\rm ker}(p)$ and the deck transformation group of this finite
cover is $G$.

Note that for any $\alpha\in G$, $\rho(\alpha)$ maps $\mbox{\rm ker}(p)$ to itself
because
\[(u,\alpha)\circ (s,I)\circ (u,\alpha)^{-1}=(\rho(\alpha)s,I)\]
and thus if $(s,I)\in\Gamma$, then $(\rho(\alpha)s,I)\in \Gamma$.

It follows that the  map $\rho\co  {\rm O}(n)\to \mbox{\rm Aut}(\re^n)$ induces a map
$\psi=co G\to \mbox{\rm Aut}(\mathbb T^n=\re^n/\mbox{\rm ker}(p))$. As an action $\phi$ of $\mathbb T^n$ on
$\re^n/\mbox{\rm ker}(p)$ we take $x\mapsto x+t$.
To see that this defines an $\Fs$--structure we check
the condition $\alpha(\phi(t)(x))=\phi(\psi(\alpha)(t))(\alpha (x))$ for $\alpha\in G$ which just says
$\alpha(x+t)=\alpha(x)+\alpha(t)$.

\subsection{Ends of hyperbolic manifolds} The following description may be found in \cite{Eb}.
Let $\Gamma\subset \is(\H^n)$ be a lattice and let $M:=\H^n/\Gamma$. If $M$ is not compact, then
it has finitely many ends (or cusps) and the ends are in one-to-one correspondence with
conjugacy classes of subgroups of $\Gamma$ that contain parabolic elements.
For each end there is a point $x\in \H(\infty)$ (the sphere at infinity) such that if we let
$\Gamma_x$ be the stabilizer of $x$, then $\Gamma_x$ consists only of parabolic elements
which leave every horosphere $L$ at $x$ invariant. The horosphere $L$ is flat with the induced metric
(this can be easily seen in the upper-half space model with horospheres given by $x_n$ constant) and thus
$N:=L/\Gamma_x$ is a compact flat manifold.
A horocyclic neighbourhood $U$ of the end is given by the projection of open horoballs in $\H^n$.
The set $U$ is a warped Riemannian product of the flat metric on $N$ and $(0,\infty)$
whose metric is given by $e^{-2t}\,ds^2_{N}+dt^2$.

\subsection{$\Fs$--structures on quotients of $\H^3\times \E$ and ends}\label{1}
Let $\Gamma$ be a lattice in $\is(\H^3\times \E)$. By the discussion in subsection \ref{lat}, there
exists $s_0$ such that $\Gamma \cap R$ contains the translations generated by $(x,t)\mapsto (x,t+s_0)$
(and only them).

Consider the projection homomorphism $\is(\H^3\times\E)\mapsto \is(\E)\mapsto \Z_2$
(recall that $\is(\E)$ is the semidirect product of $\re$ with $ {\rm O}(1)=\Z_2$).
Then we have a homomorphism $\Gamma\mapsto \Z_2$. Its kernel is an index 2 subgroup
$\Gamma_0\subset \is(\H^3)\times\re$. The manifold $M_0=\X/\Gamma_0$ is a 2-1 cover
of $M$. But $M_0$ admits a circle action since the action of $\re$, $(x,t)\mapsto (x,t+s)$
descends to a circle action on $M_0$. The action may not descend to $M$, but $M$ is still foliated by circles.
In any case we obtain in this way an $\Fs$--structure on $M$, where
$\Psi\co \Z_2\to \mbox{\rm Aut}(S^1)$ on the non trivial element of $\Z_2$ is just $t\mapsto -t$.

Let us now take a look at the ends of $M$. Let $p_1\co \is(\H^3\times \E)\to \is(\H^3)$ be the projection
on the first factor. The group $p_1(\Gamma)$ is a lattice in $\is(\H^3)$ isomorphic to
$\Gamma/\Z$. The ends of $M$ arise from the ends of the hyperbolic $3$--orbifold
$\H^3/p_{1}(\Gamma)$. Note that the action of $p_1(\Gamma)$ on $\H^3$ is not neccesarily free and the fixed
points create the orbifold nature of the quotient. By Selberg's lemma \cite{Se}  $p_1(\Gamma)$ does contain
a finite index subgroup which acts freely on $\H^3$.

For each end of $M$, there is a point $x\in \H^3(\infty)$ and a horosphere $L$
through $x$. The set $P=L\times\E$ is a copy of Euclidean $3$--space which inherits the flat metric
from $\H^3\times \E$. If we let $\Gamma_P$ be the elements of $\Gamma$ which project under $p_1$ to ${\rm Stab}(x)$, then the horocyclic neighbourhood $V=P/\Gamma_P\times (0,\infty)$ is the end of $M$.
Then, on $V$ we have a canonical $\Fs$--structure given by subsection \ref{flatend}.

Note that the $\Fs$--structure we defined on $M$ before using the $\re$--action on the $\E$--factor
is compatible with the one we just described at the ends. In fact on $V$ the $\re$--action does descend
to a circle action leaving $P/\Gamma_P$ invariant.

\subsection{Gluing $\H^3\times\E$ pieces}\label{glue}
In this subsection we suppose that $M$ is a closed orientable
geometrisable $4$--manifold with pieces modeled on $\H^3\times \E$ and we show how to put a polarised
$\Fs$--structure on $M$.

In order to prove this, the situation we need to consider is the following.
Let $M_i=\H^3\times\E/\Gamma_i$ for $i=1,2$ and suppose
$M_i$ has one end of the form $P_i\times (0,\infty)$ for $i=1,2$ and $f \co P_1\to P_2$
is a diffeomorphism. The manifolds $P_i$ are flat.
We wish to show that $M=M_1\cup_{f} M_2$ has an $\Fs$--structure.
The diffeomorphism type of $M$ only depends on the isotopy class of $f$.
We will use the fact that on a flat $3$--manifold any diffeomorphism is isotopic to an affine map, so from now
on we will suppose that $f$ is affine (this follows from either \cite{MS} or \cite{BL}).

Now according to the previous subsection we have $\Fs$--structures on each of the ends. These structures
will be compatible when the gluing map is affine. Indeed we only need to observe that
in $\re^n$, an affine map has the form $f(x)=Ax+b$, where $A$ is an invertible matrix and $b\in\re^n$
a fixed vector. Hence if we conjugate by $f$ the $\re^n$--action by translations $x\mapsto x+u$
we obtain $x\mapsto x+Au$ and these two actions commute. So we have a polarised $\Fs$--structure on $M$.

\subsection{$\Fs$--structures on quotients of $\H^2\times\E^2$ and ends.}\label{2}
Let $\Gamma$ be a lattice in $\is(\H^2\times \E^2)$. By the discussion in subsection \ref{lat} above, there
exist linearly independent vectors $w_1,w_2\in\re^2$ such that $\Gamma$ contains
the translations generated by $(x,y)\mapsto (x,y+w_i)$, for $i=1,2$ (and only them).

Consider the projection homomorphism $\is(\H^2\times\E^2)\to \is(\E^2)\to O(2)$
(recall that $\is(\E^2)$ is the semidirect product of $\re^2$ with $ {\rm O}(2)$).
Then we have a homomorphism $\Gamma\mapsto  {\rm O}(2)$ with image a finite group $G$.
Its kernel is a finite index subgroup
$\Gamma_0\subset \is(\H^2)\times\re^2$. The manifold $M_0=\X/\Gamma_0$ is a finite cover
of $M$ with $G$ as deck transformation group.
But $M_0$ admits a $2$--torus action since the action of $\re^2$, $(x,y)\mapsto (x,y+u)$
descends to a $2$--torus action on $M_0$. The action may not descend to $M$, but $M$ is still foliated by tori.
In any case we obtain in this way an $\Fs$--structure on $M$, where
$\Psi \co G\to \mbox{\rm Aut}( T^2)$ is given exactly by the $\Psi$ of subsection \ref{flatend}.

Let us now take a look at the ends of $M$. Let $p_1 \co\is(\H^2\times \E^2)\to \is(\H^2)$ be the projection
on the first factor. The group $p_1(\Gamma)$ is a lattice in $\is(\H^2)$ isomorphic to
$\Gamma/\Z^2$. The ends of $M$ arise from the ends of the hyperbolic $2$--orbifold
$\H^2/p_{1}(\Gamma)$.

For each end of $M$, there is a point $x\in \H^2(\infty)$ and a horosphere $L$
through $x$. The set $P=L\times\E^2$ is a copy of Euclidean $3$--space which inherits the flat metric
from $\H^2\times \E^2$. If we let $\Gamma_P$ be the elements of $\Gamma$ which project under $p_1$ to ${\rm Stab}(x)$, then the horocyclic neighbourhood $V=P/\Gamma_P\times (0,\infty)$ is the end of $M$.
Then, on $V$ we have a canonical $\Fs$--structure given by subsection \ref{flatend}.

Note that the $\Fs$--structure we defined on $M$ before using the $\re^2$--action on the $\E^2$--factor
is compatible with the one we just described at the ends.

\subsection{$\Fs$--structures on quotients of $\Sl\times\E$ and ends}\label{3}
Let $\Gamma$ be a lattice in $\is(\Sl\times \E)$. Since $\Sl$ does not admit orientation reversing isometries
and $M$ is orientable we see that $\Gamma\subset \is(\Sl)\times\re$.
Recall that we have the sequence
\[0\to \re\to \is(\Sl)\to\is(\H^2)\to 1\]
and $\re$ is central in $\is_{0}(\Sl)$. Hence $\is(\Sl)\times\re$ contains a copy of $\re^2$.
By the discussion in subsection \ref{lat}, there
exist linearly independent vectors $w_1,w_2\in\re^2$ such that $\Gamma$ intersects $\re^2$
in the lattice $\Z w_1+\Z w_2$.

We have a homomorphism
\[\Gamma\to \is(\Sl)\to\is(\H^2)\to \Z_2=\is(\H^2)/PSL(2,\re).\]
The kernel of this homomorphism gives an index 2 subgroup
$\Gamma_{0}\subset \is_{0}(\Sl)\times\re$
and the manifold $M_0=\X/\Gamma_0$ is a 2-1 cover of $M$.
But $M_0$ admits a $2$--torus action since the action of $\re^2$ on $\Sl\times\E$
descends to a $2$--torus action on $M_0$. The action may not descend to $M$, but $M$ is still foliated by tori.
In any case we obtain in this way an $\Fs$--structure on $M$, where
$\Psi \co \Z_2\to \mbox{\rm Aut}(\mathbb T^2)$ is given by $(t_1,t_2)\mapsto (-t_1,t_2)$.

Let us now take a look at the ends of $M$. Let $p_1 \co \is(\Sl\times \E)\to \is(\H^2)$ be the composition
of the projection on the first factor with $\is(\Sl)\to \is(\H^2)$.
The group $p_1(\Gamma)$ is a lattice in $\is(\H^2)$ isomorphic to
$\Gamma/\Z^2$. The ends of $M$ arise from the ends of the hyperbolic $2$--orbifold
$\H^2/p_{1}(\Gamma)$.

For each end of $M$, there is a point $x\in \H^2(\infty)$ and a horosphere $L$
through $x$. The Lie group $\Sl$ is an $\re$--bundle over $\H^2$, so inside $\Sl$ we now get
a copy $F$ of Euclidean $2$--space given by those $\re$--lines over $L$.

The set $P=F\times\E$ is a copy of Euclidean $3$--space which inherits the flat metric
from $\Sl\times \E$. If we let $\Gamma_P$ be the elements of $\Gamma$ which project under $p_1$ to ${\rm Stab}(x)$, then the horocyclic neighbourhood $V=P/\Gamma_P\times (0,\infty)$ is the end of $M$.
Then, on $V$ we have a canonical $\Fs$--structure given by subsection \ref{flatend}.

As in the other cases, the $\Fs$--structure we defined on $M$ before using the $\re^2$--action
on $\Sl\times\E$ is compatible with the one we just described at the ends.
In fact, on $V$ we do have a $2$--torus action leaving $P/\Gamma_P$ invariant.

\begin{Theorem}\label{mixed} A closed orientable and geometrisable 4-manifold
with pieces modelled only on $\H^3\times \E$, $\H^2\times \E^2$ or $\Sl\times\E$
admits a polarised $\Fs$--structure.
\end{Theorem}
\begin{proof}
Subsections \ref{1}, \ref{2} and \ref{3} exhibit $\Fs$--structures on each piece, such that
at the flat ends we have the canonical $\Fs$--structure defined in subsection \ref{flatend}.
Gluing by affine diffeomorphisms ensures compatibility on the overlaps as explained
in subsection \ref{glue}. By inspection we see that the structure is polarised. \end{proof}

\section{ Manifolds which decompose into $\mathbb{F}^4$--pieces}

\subsection{The geometry  $\mathbb{F}^4$}

Suppose we have  a finite volume manifold $M$ modelled on $\mathbb{F}^4$. The fundamental group $\Gamma$ of $M$ is a lattice in $\mathbb{R}^2\ltimes {\rm SL}(2,\mathbb{R})$. It must meet $\mathbb{R}^2$  in $\mathbb{Z}^2$, otherwise the volume of $M$ would not be finite. Denote by $\overline{\Gamma}$ the image of $\Gamma $ in ${\rm SL}(2,\mathbb{R})$, notice $\overline{\Gamma}= \Gamma / \mathbb{Z}^2$. We can now see that $M = \mathbb{F}^4 / \Gamma$ is an elliptic surface over $B=\mathbb{H}^2 / \overline{\Gamma}$, where $B$ is a non-compact orbifold \cite[p.150]{Wall}.

The identity component of ${\rm Iso}(\mathbb F^4)$ coincides with ${\rm Iso}^{+}(\mathbb F^4)$
and is given by the semidirect product
$\mathbb{R}^{2}\ltimes_{\alpha}{\rm SL}(2,\mathbb{R})$, with $\alpha$ the natural action of
${\rm SL}(2,\mathbb{R})$ on $\mathbb{R}^{2}$.
Let $\Gamma\subset {\rm Iso}^{+}(\mathbb F^4)$ be a lattice, so that $M=\mathbb{F}^4/\Gamma$ is a finite
volume manifold modelled on $\mathbb F^4$. Let $p \co \mathbb{R}^{2}\ltimes_{\alpha}{\rm SL}(2,\mathbb{R})\to {\rm SL}(2,\mathbb R)$
be the projection homomorphism. By a theorem of Wang (cf. \cite[8.27]{Ra}), $\Gamma$ meets $\mathbb R^2$ in a lattice
isomorphic to $\mathbb Z^2$. The quotient $\Gamma/\mathbb Z^2$ is isomorphic to $p(\Gamma)$.
As in the case of flat manifolds, the structure of semidirect product implies that
if $A\in p(\Gamma)$, then $A$ maps $\Gamma\cap \re^2$ to itself. Thus we have an induced
homomorphism $\psi \co p(\Gamma)\to {\rm Aut}(\mathbb T^2=\mathbb R^2/(\Gamma\cap\mathbb R^2))$.
The manifold $M$ is $\mathbb T^2\times\mathbb H^2$ modulo the action of $p(\Gamma)$, where it acts
on $\mathbb T^2$ via $\psi$ and on $\mathbb H^2$ in the usual way. The quotient $B:=\mathbb H^2/p(\Gamma)$
is a hyperbolic orbifold of finite volume and hence $M$ is an orbifold bundle over $B$.
If $B$ is smooth, ie $p(\Gamma)$ acts without fixed points, then $M$ is a torus bundle
over $B$ with structure group ${\rm SL}(2,\Z)$ and $\psi$ is precisely its holonomy.

\subsection{$\mathbb{F}^4$--manifolds as elliptic surfaces}

The manifold $M$ is also an elliptic surface over $B$. There are no singular fibres, but there could
be exceptional fibres with multiplicity $m\geq 2$. These would correspond to cone points on $B$ of cone
angle $2\pi/m$. The $\mathcal T$--structure is quite visible from this: $M$ is obtained from a torus bundle
after perhaps some logarithmic transforms.

\subsection{Ends of $\mathbb{F}^4$--manifolds}

The ends of $M$ arise from parabolics elements in $p(\Gamma)$. If $L\subset \mathbb H^2$ is an appropriate horosphere
left invariant by a parabolic element $A\in p(\Gamma)$, then the cusp will have the form
$P\times (0,\infty)$, where $P=(\mathbb T^2\times L)/\Z$, where $\Z$ is generated by $A$.
This exhibits the boundary of the ends as torus bundles over the circle.

\subsection{Affine transformations of Lie groups}Let $G$ be a Lie group and ${\rm Aut}(G)$ be the group of continuous automorphisms of $G$. Then the group ${\rm Aff}(G)$ of affine transformations of $G$ is isomorphic to the semi-direct product $A(G):=G\ltimes {\rm Aut}(G)$ with the operation,
 $$(g_1, \alpha_1)(g_2, \alpha_2)=(g_1 \alpha_1(g_2), \alpha_1 \alpha_2), \quad g_1,g_2 \in G, \quad \alpha_{i}\in {\rm Aut}(G).$$

It has a Lie group structure and acts on $G$ by $(g, \alpha)x = g\alpha(x)$ for $(g, \alpha) \in A(G)$, $x\in G$.

The left inverse of $(g, \alpha)$ is $(g, \alpha)^{-1} = ((\alpha^{-1}(g))^{-1}, \alpha^{-1})$.
\begin{eqnarray*}
(g, \alpha)^{-1}(g, \alpha) &=& ((\alpha^{-1}(g))^{-1}, \alpha^{-1})(g, \alpha)  \\
                                         &=& ((\alpha^{-1}(g))^{-1} \alpha^{-1}(g), \alpha^{-1}\alpha ) \\
                                         &=& (e, Id ).
\end{eqnarray*}

It was first noticed by Kamber--Tondeur in \cite{KT} that the action of $A(G)$ on $G$ defines an isomorphism $i \co A(G)\rightarrow {\rm Aff}(G)$.
%
%
%
 The following lemma is useful for computations.
 \begin{Lemma}\label{nilaff}  {Let $G$ be a Lie group, $\rho$ and $\sigma $ elements of the centre of $L$ and $A\in {\rm Aff}(G)$, then $A^{-1}\rho A \sigma = \sigma A^{-1}\rho A$.}
 \end{Lemma}
 \begin{proof} Let $\rho_{A}= A^{-1}\rho A$, where $A=(g, \alpha)$ in $A(G) \cong {\rm Aff}(G)$ and similarly $\rho=(\rho, Id)$ and $\sigma=(\sigma, Id)$.

 If $\rho $ is in the centre of $G$ then for an $\alpha$ in ${\rm Aut }(G)$ and a $g$ in $G$ we have that
  \[ \rho g = g \rho  \Rightarrow \alpha(\rho g)= \alpha(g \rho ) \Rightarrow \alpha(\rho) \alpha(g) = \alpha (g) \alpha(\rho ).\]

 We now compose the above elements to see that $\rho_{A} =(\alpha^{-1}(\rho), Id):$
 \begin{eqnarray*}
\rho_{A} &=& A^{-1}\rho A \\
       &=& A^{-1}[(\rho , Id)\circ ( g , \alpha)] \\
       &=&   A^{-1}\circ (\rho g, \alpha) \\
       &=& ((\alpha^{-1}(g))^{-1}, \alpha^{-1})\circ (\rho g, \alpha) \\
       &=& ((\alpha^{-1}(g))^{-1}\alpha^{-1}(\rho g), \alpha^{-1} \alpha )\\
       &=& ((\alpha^{-1}(g))^{-1}\alpha^{-1}(g \rho), Id))\\
       &=& ((\alpha^{-1}(g))^{-1}\alpha^{-1}(g)\alpha^{-1}( \rho), Id) )\\
       &=& (\alpha^{-1}(\rho), Id)
\end{eqnarray*}
The above calculation implies:
\begin{eqnarray*}
\rho_{A}\sigma &=&  (\alpha^{-1}(\rho), Id)\circ (\sigma, Id) \\
                         &=& (\alpha^{-1}(\rho) \sigma, Id) \\
                         &=& (\sigma \alpha^{-1}(\rho), Id) \\
                         &=& (\sigma, Id)\circ (\alpha^{-1}(\rho), Id)\\
                         &=& \sigma \rho_{A}
\end{eqnarray*}
Therefore $A^{-1}\rho A \sigma = \sigma A^{-1}\rho A.$
\end{proof}

\subsection{$\Ts$--structures on manifolds with decomposition into $\mathbb{F}^4$--pieces}

\begin{Theorem}\label{f4} Let $M$ be a closed orientable complete manifold with a geometric decomposition into orientable pieces modelled only on $\mathbb{F}^4$, then $M$ admits a $\mathcal{T}$--structure.
\end{Theorem}

\begin{proof}  First we will see how the $\mathbb{F}^4$--pieces  of the geometric decomposition on $M$ admit a $\mathcal{T}$--structure.

Let $N$ denote one such $\mathbb{F}^4$--piece, then $N$ is an open elliptic surface over the base $B$. Let $m$ be the number of cusps of $B$ and $p_{i}$ one such cusp.  Denote by $E_{i}$ the end of $N$ corresponding to the cusp $p_{i}$ of $B$. We know that $E_{i}$ is a $Nil^3$--manifold and $\pi_1(E)$ is isomorphic to $\Gamma_{k_{i}}$ as above, for some $k_{i}\in \Z$.

Consider a small horocyclic neighbourhood $U$ of $p_{i}$, let $B^0:= B -U$ and $N^0 \rightarrow B^0$ be the corresponding elliptic surface obtained by restriction. Identify the boundary $\partial N^0$ of $N^0$ with itself using the identity to form the double $DN^0$ of $N^0$.

Now $DN^0$ is a compact elliptic surface over the double of $B^0$, so $DN^0$ admits a $\mathcal{T}$--structure whose orbits are the elliptic fibres \cite[Thm.5.10]{PP}.

When we restrict the $\mathcal{T}$-structure on $DN^0$ to $N^0$ we obtain a $\mathcal{T}$--structure on $N^0$.

 Recall that  $E_{i}$ is a $T^2$--bundle over $S^1$ with geometric monodromy  $$\mu_{i} = \left(\begin{array}{cc}
1 & k_{i} \\
0 & 1  \\  \end{array} \right) \in {\rm SL}(2, \mathbb{R}).$$ The monodromy around $\partial U$, the boundary of $U$, is also $\mu_{i}$. This allows us to extend the
$\mathcal{T}$--structure on $N^0$ to $N$, because the action of $T^2$ on the elliptic fibres behaves well with respect to any element of ${\rm SL}(2, \Z)$. This holds in particular for $\mu_{i}$, the corresponding actions will commute after conjugation by  $\mu_{i}$.

A collar nieghbourhood of $E_{i}$ in $N$ is diffeomorphic to $E_{i}\times (0, \infty)$, which is $Nil^3 / \Gamma_{k_{i}}\times (0, \infty) : = V$.

Both actions leave $V$ invariant, as they leave every slice $E_{i}\times \{ t \} $ invariant for every $t \in (0, \infty ) $.


As translations along the $z$ axis in $Nil^3$, given by $[x,y,z]\stackrel{\sigma}{\mapsto} [x,y,z + \sigma]$, $\sigma \in \re$ are central in $Nil^3$, they descend to an $S^1$--action $\sigma$ on $V$,
$$([x,y,z],t)\stackrel{\sigma}{\mapsto} ([x,y,z + \sigma],t).$$

On $V$ the $\mathcal{T}$--structure $\tau$ defined in $N$ by rotation on the elliptic fibres commutes with $\sigma$ on $N\cap V$,
\begin{eqnarray*}
\tau \sigma ([x,y,z],t) &=& \tau ([x,y,z + \sigma],t) \\
                                 &=& ([x + \tau _1 , y +\tau_2 ,z + \sigma],t) \\
                                 &=& \sigma ([x + \tau _1 , y +\tau_2 ,z ],t) \\
                                 &=& \sigma \tau ([x,y,z],t).
\end{eqnarray*}

 Assume $N_1$ and $N_2$ are two $\mathbb{F}^4$--manifolds which are glued along  $E_1$ and $E_2$, components of their respective boundaries.

 Let $h \co E_1\rightarrow E_2$ be the glueing diffeomorphism; we will see in the next section that $h$ is isotopic to an affine transformation $\alpha \co E_1\rightarrow E_2$.

 When we use isotopic diffeomorphisms to identify boundary components, we obtain diffeomorphic manifolds  \cite{Hir}. So it is enough to work with $\alpha$, as we are interested in the existence of a $\mathcal{T}$--structure up to diffeomorphism.

 Define $\rho_{\alpha}:=\alpha^{-1}\rho \alpha$, we need to show that $  \rho_{\alpha}\sigma= \sigma \rho_{\alpha}$ for the $S^1$--actions $\rho , \sigma $ on $E_1$ and $E_2$ which are induced by translations along the $z$--axis of $Nil^3$.

 The affine transformation $\alpha$ lifts to an affine transformation $A$ of $Nil^3$ which is $\pi_1(E_1)$--invariant, that is $A$ sends $\pi_1(E_1)$--orbits to  $\pi_1(E_1)$--orbits in $Nil^3$.

 Both $\rho$ and $\sigma$ lift to translations along the $z$--axis of $Nil^3$, which we will call $\tilde{\rho}$ and $\tilde{\sigma}$.

 Therefore $\rho_{\alpha}\sigma= \sigma \rho_{\alpha}$ follows from $\tilde{\rho_{A}}\tilde{\sigma}= \tilde{\sigma} \tilde{\rho_{A}}$ on $Nil^3$, where again $\tilde{\rho_{A}}:=A^{-1}\tilde{\rho}A$. This was shown in lemma \ref{nilaff}, as both  $\tilde{\rho}$ and $\tilde{\sigma}$ are central in $Nil^3$.

 By repeating the same procedure on each geometric piece of $M$ and on each pair of identified cusps, we give $M$ a   $\mathcal{T}$-structure.
 \end{proof}

\section{Diffeomorphisms of closed flat $3$--manifolds and $Nil^3$--manifolds}\label{affdiff}

In his review on problems in low dimensional topology \cite[p.137]{Kir}  R Kirby points out that the following is a consequence of the results of P Scott and W Meeks \cite{MS}. Let $M$ be a $3$--manifold modelled on $\re^3$ or $ Nil^{3}$, ${\rm Aff}(M)$ denote the group of affine transformations of $M$ and ${\rm Diff}(M)$ the group of diffeomorphisms of $M$. The inclusion ${\rm Aff}(M) \hookrightarrow {\rm Diff}(M)$  induces an isomorphism on components
\[ \pi_0( {\rm Aff}(M))\stackrel{\cong}{ \hookrightarrow} \pi_0( {\rm Diff}(M)).\]
\begin{rem}
It follows that every diffeomorphism of a compact  flat  $3$--manifold or $Nil^{3}$--manifold is isotopic to an affine transformation.
\end{rem}

 A proof of this result can be reproduced by induction from the results on periodic diffeomorphisms shown in  \cite{MS}---this was communicated to the author by Scott. An alternative approach has been suggested by A. Verjovsky \cite{Ver},  which uses the following deep fact about closed one forms on 3-manifolds. Laudenbach--Blank showed in \cite{BL}that two closed non-singular $1$--forms on a $3$--manifold are isotopic if and only if they are cohomologous.

\begin{Corollary}\label{3diff}  Any diffeomorphism $g$ of $T^3$ is isotopic to an affine transformation.
\end{Corollary}
\begin{proof} {\rm{ (Verjovsky)}} The map $g$ induces a linear map on $H^1(T^3)$, composing $g$ with the inverse of this linear map, we can assume that $g^{\ast}\co H^1(T^3)\rightarrow H^1(T^3)$ is the identity. Let $p\co  T^3=S^1\times S^1 \times S^1 \rightarrow S^1$ be the projection onto the first factor and let $\omega = \pi^{\ast }(d\theta)$, where $d\theta$ denotes the metric on $S^1$. The form $\omega $ is a closed non-singular $1$--form on $T^3$. Since $g^{\ast }\omega$ is cohomologous to $\omega$, by the Laudenbach--Blank theorem, $g^{\ast }\omega$ and $\omega$ are isotopic through an isotopy $(h_{t})_{\{ 0 \leq t \leq 1 \} }$ such that $h_{0}= \rm{Id}$ and $h_{1}^{\ast}(g^{\ast}\omega)=\omega$.

The map $f=g\circ h_{1}$, fixes  $\omega $ in cohomology and therefore fixes each torus $T_{\theta}= \{ e^{2\pi i \theta} \} \times T^2$. Let $f_{\theta} $ be the restriction of $f$ to $T_{\theta}$. The map $\theta \mapsto f_{\theta}$ defines a loop $S^1\rightarrow {\rm{Diff}}(T^2)$; in fact the image of this loop lies in ${\rm{Diff_{0}}}(T^2)$, the subgroup of diffeomorphisms isotopic to the identity, because $f_{\theta}$ induces the identity in the cohomology of $T^2$. Because ${\rm{Diff_{0}}}(T^2)$ retracts to the group of translations and therefore we can retract out loop to a map $S^1 \rightarrow T^2$. This map is homotopic to a constant, since it is the identity in (co)-homology. \end{proof}
The same result holds for compact $\re^3$ and $ Nil^{3}$ manifolds, following analogous arguments.

\section{Fundamental groups of geometrisable manifolds}

The following topological fact allows us to describe the fundamental group of a manifold with a proper geometric decomposition \cite[p.35]{Sta}.

Assume $K_{A}$ and $K_{B}$ are $n$-dimensional topological complexes having fundamental groups $A$ and $B$ and intersecting in a connected subcomplex having fundamental group $C=A \cap B$ such that $C$ injects into both $A$ and $B$. The Seifert-van Kampen Theorem \cite{Sei2} \& \cite{Van} states that the fundamental group of $K_{A} \cup K_{B}$ is $A\ast_{C} B$.

If $K_{A}$  contains two isomorphic subcomplexes $K_{C}$ and $K_{C'}$ with fundamental groups $C$ and $\phi(C)$ (where $\phi$ is the map on $\pi_1$ by the isomorphism between $K_{C}$ and $K_{C'}$) and we attach to $K_{A}$ the cylinder $K_{C}\times [0,1]$ identifying $K_{C}$ with $K_{C}\times 0$ and $K_{C'}$ with $K_{C'}\times 1$. It also follows form the Seifert-van Kampen Theorem that we obtain a complex with fundamental group $A\ast_{C}^{\phi}$. The notation  $A\ast_{C}^{\phi}$ denotes the amalgamation of $A$ over the image of $\phi$.

Let $M$ be an orientable smooth four manifold which admits a proper geometric decomposition. It follows from the discussion above  that $\pi_1(M)$ is isomorphic to an amalgamated product $A\ast_{C}B$ or to an  extension $A\ast_{C}^{\phi}$, where $A$ is the fundamental group of one of the geometric pieces.

A free product with amalgamation $A\ast_{C}B$ where $C$ is a subgroup of both $A$ and $B$ is \emph{non-dihedral} if the two inclusions $C\subset A$ and $C \subset B$ are strict and if, moreover, the index of $C$ is not 2 in both $A$ and $B$. An extension  $A\ast_{C}^{\phi}$, where $\phi$ is an isomorphism from some subgroup $C$ of $A$ onto some subgroup $C'$ of $A$ is \emph{non-semi direct} if at least one of the inclusions $C\subset A$ or $C' \subset A$ is strict.

It was shown  by P. De la Harpe that if a group $\Gamma $ is isomorphic to either a non-dihedral amalgamated product $A\ast_{C}B$ or to a non-semi direct extension $A\ast_{C}^{\phi}$, then $\Gamma $ is of exponential growth \cite{DlH}.

\begin{Lemma}\label{d-eg} {The fundamental group of a smooth four manifold $M$ with a proper geometric decomposition  has exponential growth.}
\end{Lemma}
\begin{proof} Follows from the previous discussion. \end{proof}

Therefore if $M$ admits a proper geometric decomposition, for any smooth Riemannian metric $g$ on $M$ we have that $\h >0$.

Recall that the fundamental group of a connected sum is the free product of the fundamental groups of the summands.

If $A$ and $B$ are two finitely generated groups, then the free product $A\ast B$ contains a free product of rank $2$ unless $A$ is trivial or $B$ is trivial or $A$ and $B$ both have order $2$. Therefore if $M$ and $N$ are differentiable manifolds with $\pi_1(M)=A$ and $\pi_1(N)=B$, then $\pi_1(M\# N)$ will grow exponentially and again $\h >0$ for any smooth metric $g$ on $M\# N$, unless $\pi_1(M)$ and $\pi_1(N)$ are trivial, $\pi_1(M)$ is trivial and $\pi_1(N)$ has order $2$, or $\pi_1(M)$ and $\pi_1(N)$ both have order $2$.

\section{ Proof of the Main Results}

\subsection{Proof of Theorem 1}

\begin{proof}
 In each case we can construct an $\Fs$--structure.

    For $\mathbb{S}^{4}$ and $\mathbb{C}{\rm P}^{2}$, the only manifolds
    modelled on these geometries are $S^{4}$ and $\mathbb{C}{\rm P}^{2}$.
    These two manifolds have $S^{1}$--actions so they admit a
     $\Ts$--structure.

    If $M$ is modelled on $\mathbb{S}^{3}\times \mathbb{E} , \mathbb{H}^{3}\times \mathbb{E} ,
\Sl\times \mathbb{E} , \mathbb{N}il^{3}\times \mathbb{E} ,
\mathbb{N}il^{4}$ or $ \mathbb{S}ol^{4}_{1}$, then $M$ is foliated by geodesic circles. By  Theorem \ref{s1fol}, this foliation allows us to define a  polarised $\Ts$--structure.

In the case of $M$ being modelled on $\mathbb{S}^{2} \times
\mathbb{E}^{2}$ or $\mathbb{H}^{2}\times\mathbb{E}^{2}$, then $M$ is
Seifert fibred. Hence $M$ admits a polarised $\Ts$--structure. We saw in Theorem \ref{seifert} how to define circle actions in the fibres which behave well with respect to the Seifert fibration, in that they commute with the structure group and so define a $\Ts$--structure.

When $M$ is modelled on $\mathbb{S}ol^{4}_{m,n}$ or
$\mathbb{S}ol^{4}_{0}$, $M$ is actually diffeomorphic to a mapping torus of
$T^{3}$. With this description $M$ can be given a polarised $\Ts$--structure, as was explained in Theorem \ref{maptorus}.

 C\,T\,C Wall has shown that closed orientable $\mathbb{S}^{2}\times
\mathbb{S}^{2}$ and  $\mathbb{S}^{2}\times
\mathbb{H}^{2}$ manifolds $M$ are diffeomorphic to  complex ruled surfaces \cite{Wall}. We have constructed in Theorem \ref{s2foliations} a $\mathcal{T}$--structure on foliated manifolds whose leaves are $S^2$ or $\mathbb{R}P^2$, which include all these cases. The idea here is that the $S^2$ leaves can be rotated consistently, endowing $M$ with a  $\mathcal{T}$--structure. In general this $\Ts$--structure will not be polarised because for some of these manifolds $\chi (M)>0$ \cite{Hil} and the vanishing of the minimal volume implies that all the characteristic numbers of $M$ are zero.

Finally, the flat case, there are 27 orientable and 48 non-orientable compact flat 4-manifolds, out of which only 8 admit a complex structure \cite[p.199]{BHPV} and with 3 exceptions they are all Seifert fibred \cite[p.146]{Hil}. As we mentioned in section \ref{flat} one of Bieberbach's theorems implies that all flat manifolds admit an $\Fs$--structure because a closed flat manifold is finitely covered by a torus of the same dimension.

Therefore if $M$ is modelled on a geometry in {\bf V} then $M$ admits a $\Ts$--structure.
\end{proof}

\subsection{Proof of Theorem A}

\begin{proof}    $[$ $i) \Rightarrow ii)$ $]$ This is content of Theorem 1.

[ $ii) \Rightarrow iii)$ ] This is Theorem A of Paternain and Petean in \cite{PP}.

[ $iii) \Rightarrow iv)$ ] Follows directly from the string of inequalities between the asymptotic invariants mentioned in the introduction.


[ $iv) \Rightarrow v)$ ] If $M$ is modelled on ${\bf H}$ then $||M||>0$ by Proposition \ref{Hpos}. Therefore $||M||=0$ implies $M$ is modelled on a geometry in ${\bf V}$.

It follows from Theorem 1 that $M$ admits a $\Ts$--structure and by Paternain and Petean's Theorem A in \cite{PP} $M$ collapses with curvature bounded from below.

[ $v) \Rightarrow i)$ ] Again we will show the contrary. Let $M$ be a manifold modelled on a geometry in ${\bf H}$, then Proposition \ref{Hpos} implies $||M||>0$. As a consequence of the results in \cite{PP} and \cite{CG}, $||M||$ bounds ${\rm Vol_{K}} (M)$ (up to some constants depending only on the dimension $n$) from below. Therefore for some constant $c_{n}$ we have
\[ 0<||M||\leq c_{n} {\rm Vol_{K}} (M) \] so $M$ can not collapse with curvature bounded from below.
 \end{proof}

 \subsection{Proof of Theorem B}

\begin{proof}If the connected sum components of $M$ admit a $\Ts$-structure then this extends to $M$ under the connected sum. The same is true for $\Fs$-structures if one of the open sets of the $\Fs$-structure has a trivial covering \cite[p.437]{PP}. In all the cases in ${\bf V}$ we may achieve this as can be seen from the proofs of theorems 1,  \ref{s2foliations}, \ref{mixed} and \ref{f4}. \end{proof}

%
%
%
%

\end{document}